\documentclass[12pt]{article}

\usepackage{amsbsy,amssymb,amscd,amsfonts,latexsym,amstext,delarray,
amsmath,amsgen, amsthm, epsfig}
 
\usepackage{epsfig}
\input psfig.sty

\newtheorem{theorem}{Theorem}

\newtheorem{proposition}[theorem]{Proposition}
\newtheorem{lemma}[theorem]{Lemma}

\newtheorem{definition}[theorem]{Definition}

\newtheorem{remark}[theorem]{Remark}

 \newcommand{\ie}{{\it i.e.\/}\ } 

\newcommand{\ra}{\rightarrow}
\newcommand{\fl}{\forall}
\newcommand{\wt}{\widetilde}

\newcommand{\s}{\sigma}
\newcommand{\D}{\Delta}

\newcommand{\Zb}{\mathbb{Z}}

\newcommand{\ot}{\otimes}

\newcommand{\Hc}{\mathcal{H}}
\newcommand{\g}{\gamma}
\newcommand{\vp}{\varphi}
\newcommand{\ve}{\varepsilon}
\newcommand{\Cb}{\mathbb{C}}

\newcommand{\FA}{\mathfrak{A}}

\newcommand{\FH}{{\mathfrak{H}}}

\newcommand{\Fa}{\mathfrak{a}}
\newcommand{\FG}{{\mathfrak{g}}}

\newcommand{\Fg}{\mathfrak{g}}
\newcommand{\Fh}{\mathfrak{h}}

\newcommand{\Fl}{\mathfrak{l}}
\newcommand{\Fo}{\mathfrak{o}}

\def\Qb{{\mathbb Q}}
\def\Cb{{\mathbb C}}

\def\Nb{{\mathbb N}}
\def\Rb{{\mathbb R}}

\def\Zb{{\mathbb Z}}

\def\Ac{{\mathcal A}}

\def\Cc{{\mathcal C}}

\def\Fc{{\mathcal F}}
\def\Gc{{\mathcal G}}
\def\Hc{{\mathcal H}}

\def\Kc{{\mathcal K}}
\def\Lc{{\mathcal L}}

\def\Nc{{\mathcal N}}

\def\Vc{{\mathcal V}}

\def\a{\alpha}
\def\b{\beta}

\def\d{\delta}

\def\g{\gamma}
\def\k{\kappa}

\def\lb{\lambda}
\def\om{\omega}
\def\s{\sigma}
\def\t{\theta}
\def\ve{\varepsilon}

\def\vp{\varphi}

\def\D{\Delta}
\def\G{\Gamma}
\def\Lb{\Lambda}
\def\Om{\Omega}

\def\fl{\forall}
\def\ify{\infty}

\def\nb{\nabla}
\def\op{\oplus}
\def\ot{\otimes}
\def\part{\partial}

\def\sbs{\subset}
\def\semi{\rtimes}

\def\ts{\times}
\def\wdg{\wedge}

\def\ra{\rightarrow}
\def\longra{\longrightarrow}

\def\text{\hbox}

\def\fl{\forall}
\def\ify{\infty}
\def\mpo{\mapsto}
\def\nb{\nabla}
\def\ot{\otimes}

\def\ra{\rightarrow}

\def\sbs{\subset}

\def\ts{\times}
\def\wdg{\wedge}
\def\wt{\widetilde}

\def\Ad{\mathop{\rm Ad}\nolimits}

\def\Diff{\mathop{\rm Diff}\nolimits}

\def\Ext{\mathop{\rm Ext}\nolimits}
\def\Hom{\mathop{\rm Hom}\nolimits}
\def\GV{\mathop{\rm GV}\nolimits}

\def\SO{\mathop{\rm SO}\nolimits}
\def\Id{\mathop{\rm Id}\nolimits}
\def\exp{\mathop{\rm exp}\nolimits}
\def\mod{\mathop{\rm mod}\nolimits}
\def\pt{\mathop{\rm pt}\nolimits}
\def\Index{\mathop{\rm Index}\nolimits}

\def\Ker{\mathop{\rm Ker}\nolimits}

\def\sign{\mathop{\rm sign}\nolimits}

\def\build#1_#2^#3{\mathrel{
\mathop{\kern 0pt#1}\limits_{#2}^{#3}}}

\def\limproj{\mathop{\oalign{lim\cr
\hidewidth$\longleftarrow$\hidewidth\cr}}}

\numberwithin{equation}{section}
\begin{document}

\title{{\sc Background independent geometry and 
 Hopf cyclic cohomology}}

\author{Alain Connes \\
        Coll\`ege de France \\
        3 rue d'Ulm \\
        75005 Paris, France \\
\and
        Henri Moscovici\thanks{Research
    supported by the National Science Foundation
    award no. DMS-0245481.} \\
    Department of Mathematics \\
    The Ohio State University \\
    Columbus, OH 43210, USA
    }

\date{ \ }
       
\maketitle

\begin{abstract}

\noindent
This is primarily a survey of the way in which Hopf cyclic cohomology  has
emerged and evolved, in close relationship with the application of the noncommutative
local index formula to transverse index theory on foliations. Being  Diff-invariant, the
geometric framework that allowed us to treat the `space of leaves' of a general
foliation provides a
`background independent' set-up for geometry that could be of relevance
to the handling of the the background independence problem in
quantum gravity. With this potential association in mind,
we have added some new
material, which complements the original paper and is 
also meant to facilitate its understanding. Section 2 gives a detailed
description of the Hopf algebra that controls the `affine' transverse geometry
of codimension $n$ foliations, and Section 5 treats 
the relative version of Hopf cyclic cohomology  in full generality, including the case of 
Hopf pairs with noncompact isotropy.

\end{abstract}

\newpage

\section*{Introduction}

Coincidentally, or perhaps as a reflection of ``ontogeny recapitulates phylogeny"
in the world of mathematical ideas, 
the brand of cyclic cohomology related to Hopf algebras came
about in a strikingly similar fashion to cyclic cohomology itself, both being
motivated by the index theory of abstract elliptic operators, at successive
stages.
\medskip

The original impetus for the development of
cyclic cohomology (as enunciated in \cite{Cobw}, see also \cite{Khal},
and presented in \cite{Cndg}), 
was to construct invariants 
for $K$-theory classes that perform the function of the classical Chern-Weil
theory in the general framework of operator algebras. Starting from
the index pairing between the $K$-homology class of a $p$-summable
Fredholm module $({\FH}, \g, F)$ and
the $K$-theory class $[e] \in K_{0}(\Ac)$ of
an idempotent in an involutive algebra $\Ac \sbs \Lc ({\FH})$,
in the graded case to fix the ideas,
\begin{equation} \label{ip}
\Index (e\, F^{+} \, e)  \, = \, (-1)^n \,Tr(\g \, e [F, e]^{2n}) \, , \quad \fl \, 2n \geq p \geq 1 \, ,
\end{equation}
one arrived, by regarding $ de = i [F, e]$ as a `quantized' differential,
to the multilinearized form of the right hand side,
\begin{equation} \label{gc}
\tau_F (a^0 ,a^1 ,\ldots ,a^{2n}) =
       Tr \left( \g a^0 [F,a^1]
\ldots [F,a^{2n}] \right), \qquad a^i \in {\cal A},
\end{equation}
that turned out to encode the quintessential features 
of the cyclic cohomology theory
for algebras.
\medskip

The  non-additive category of algebras and
algebra homomorphisms was replaced in  \cite{Cext} 
by  the additive category of modules over the cyclic category $\Lb$,
allowing the realization of cyclic cohomology as an $\Ext$ functor.
This enlargement of the scope of the theory played an essential role
years later, when the authors were faced with the formidable looking
task of concretely computing in the geometrically interesting case of foliations
the `theoretical' answer provided by the universal local index formula \cite{CMlif}.     
The gist of that formula is that, in the unbounded version of the index pairing
(\ref{ip}), it replaces the `global' cocycle (\ref{gc}) by a universal finite linear
combination of `local' cocycles of the form
\begin{equation} \label{lc}
\phi (a^0 , \ldots , a^m) = \, \int \!\!\!\!\!\! - \, a^0 [D ,a^1]^{(k_1)}
\ldots
[D , a^m]^{(k_m)} \, |D|^{- (m + 2 \vert {\bf k} \vert)} \, ,  
\end{equation}
where $T^{(k)}$ stands for the $k$th iterated commutator
of the operator $T$ with $D^2$ and 
$\displaystyle \int \!\!\!\!\!\! - \,$ is an extension of the Dixmier 
trace given by residues of spectral zeta-functions. 
\medskip

In the case
of transversely hypoelliptic  operators on foliations, algebraic
manipulations with the 
commutators appearing in (\ref{lc}) led to the emergence 
of the Hopf algebra $\Hc_n$, that plays for the transverse frame bundle
to a foliation the role of the affine group
of the frame bundle 
to a manifold. 
\medskip

Recognizing the cyclic module structure associated
to the Hopf algebra $\Hc_n$, and intrinsically related to the
`characteristic' cochains (\ref{lc}), provided 
precisely the missing principle to organize the computation.
We settled the index problem in \cite{CMhti}, as briefly sketched in $\S 1$,
by proving that the cyclic cohomology of the above cyclic module
is in fact isomorphic to the Gelfand-Fuks cohomology, in both the
`absolute' and the  `relative' case. This isomorphism is  
concretely illustrated in the codimension $1$ case in  $\S 6$, for the
Godbillon-Vey class and also for the transverse fundamental class. 
\medskip

The emergent Hopf cyclic structure applies to arbitrary Hopf algebras,
 in particular to quantum groups,
 and gives rise to characteristic classes associated to Hopf 
 actions, cf. \cite{CMhti, CMcchs}, also $\S 3$ and $\S 4$ below.
 The algebraic machinery developed in the process  
has been extended by Hajac-Kahalkhali-Rangipour-Sommerh\"auser 
to a theory with coefficients \cite{HKRS1, HKRS2}.
The characteristic map associated to a Hopf module algebra with
invariant trace has been generalized to the case of
higher traces by
Crainic \cite{Cr} and by Gorokhovsky \cite{Gor}.
It was further extended by Khalkhali and Rangipour \cite{KR},
who upgraded it to cup products in Hopf-cyclic cohomology.
For these developments we refer the reader to the cited papers.
\medskip

The geometric framework that allowed us to treat the `space of leaves' of a general
foliation is Diff-invariant and therefore provides a
`background independent' set-up for geometry that could be of relevance
in dealing with the background independence problem in
quantum gravity. With this potential association in mind,
we have added some new
material. In $\S 2$ we give a detailed
description of the Hopf algebra $\Hc_n$ and of its `standard'
module-algebra representation, while
$\S 5$ treats
the relative version of Hopf cyclic cohomology in full generality; thus,
besides the relative Hopf cyclic cohomology of the pair $(\Hc_n, \Fo_n)$, 
that has played a crucial role in understanding the Chern character of the hypoelliptic
signature operator, one can now handle pairs with noncompact isotropy, such as
 $(\Hc_{n+1}, \Fo_{n, 1})$.
\bigskip

\tableofcontents

\bigskip

\section{Background independent geometry and the local index formula}

In the noncommutative approach
a  geometric structure on a `space' is specified
by means of a spectral triple $(\Ac ,{\FH} ,D)$.
$\Ac$ is an involutive algebra of bounded operators in a Hilbert space
${\FH}$, and represents the `local coordinates' of the space.
 $D$ is an unbounded selfadjoint operator on ${\FH}$, which
 has bounded commutators with the `coordinates', and whose
 inverse $D^{-1}$
corresponds to the infinitesimal line  element  $ds$
in Riemannian geometry.  In addition to
its metric significance,  $D$ carries an important 
topological meaning, that of a $K$-homology cycle which
 represents the fundamental class
of the `space' which is the spectrum of $\Ac$.

When this space is an ordinary spin manifold $M$, 
\ie when the algebra $\Ac=C^\infty(M)$, one obtains
a natural spectral triple $(\Ac ,{\FH} ,D)$ by
fixing a Riemannian metric on $M$
and taking for  $D$ the Dirac operator 
in the Hilbert space ${\FH}$ of square integrable spinors.
\medskip

At first sight,  \ie when viewed from the classical
viewpoint of Riemannian geometry,
the transition from the local differential geometric
set-up to the 
operator theoretic framework might appear as a mere
translation. That this is far from being the case,
even in the classical framework, requires an explanation
which we now give below.

In order to define the transverse geometry,
\ie the geometry of the `space' of leaves,
for a general foliation, one is confronted with the 
problem of finding a geometric structure that is 
\textit{invariant} under all diffeomorphisms of a
given manifold $M$. Indeed, the action of the 
holonomy on a complete transversal $M$ to a foliation 
is  as wild (in general) as that of an arbitrary
(countable) subgroup of $\Diff(M)$,
and invariance under holonomy
is a necessary constraint 
when passing to the space of leaves.

The standard geometric notions are of course 
\textit{equivariant} with respect to $\Diff(M)$,
but they are not \textit{invariant} . 
In fact, it is well known that the 
group of isometries of a Riemannian manifold $N$
is a finite dimensional Lie group and is 
thus incomparably smaller than the group
$\Diff(M)$ of any manifold.
\medskip

The first virtue of the operator theoretic framework of
noncommutative geometry is that it only requires
invariance to
hold at the level  of the \textit{principal symbol}
(in classical pseudodifferential
terms) of the operator $D$.
When $D$ is an elliptic operator the gain is non-existent
since in that case the symbol specifies the metric.
But the first main point is that the theory applies 
with no change when $D$ is only \textit{hypoelliptic},
and this allows to treat `para-Riemannian'
spaces, which admit groups of isometries as large as
diffeomorphism groups.

This allows to handle 
 $n$-dimensional geometry in the following ``background independent"
way \cite{CMlif}. One first replaces a given manifold $M^n$ (with no extra 
structure except an orientation) by the total
space of the bundle $PM = F^{+}M /SO (n)$, where $F^{+}M$ is 
the $GL^{+}(n, \Rb)$--principal bundle of oriented frames on
$M^n$. The sections of $\pi: PM  \ra M $ are precisely the
Riemannian metrics on $M$ but unlike the 
space of such metrics the space $P$ is still
a finite dimensional manifold. 
The total space $PM$ itself admits a canonical, and thus
$\Diff^{+}(M)$-invariant, `para-Riemannian'
structure, which can be described as follows.
The vertical subbundle $\Vc \sbs T(PM)$, $\Vc =\Ker \pi_*$, 
carries natural Euclidean structures on each of its
fibers, determined solely by 
fixing once and for all a choice of a $GL^+(n,  \Rb)$-invariant
Riemannian metric on 
the symmetric space $GL^{+}(n,  \Rb) / SO (n)$. On the other hand,
the quotient bundle $\Nc = T(PM)/ \Vc$ comes equipped with a
tautologically defined Riemannian structure: every point $p\in PM$ is an
Euclidean structure on $T_{\pi (p)} (M)$ which is identified to 
$\Nc_{p}$ via $\pi_*$. 

Since no non-canonical choice were involved so far, the obtained
structure on $PM $ is \textit{invariant} under the canonical
lift of the action of $\Diff^{+}(M)$. In particular
any hypoelliptic operator whose  \textit{principal symbol}
only depends upon the above `para-Riemannian' structure
will have the required invariance to yield a 
spectral triple governing the geometry
in a ``background independent"
manner. Since the object of our interest is the $K$-homology class
of the spectral triple (and we can freely use the Thom
isomorphism to pass from the base $M$ to the 
total space $PM$ in an invariant manner), we shall
take for $D$ the hypoelliptic signature operator.
The precise construction of $D$, to be recalled below,
involves the choice of a connection on the frame bundle but this 
choice does not affect the  \textit{principal symbol}
of $D$ and thus plays an innocent role which does not
alter the fundamental $\Diff^{+}(M)$-invariance
of the spectral triple. More precisely, we have shown
in \cite{CMlif} that it does define in full 
generality  a spectral triple on the crossed 
product of $PM$ by $\Diff^{+}(M)$.
 \smallskip

It is worth mentioning at this point that this
construction, besides allowing to handle arbitrary foliations,
could be of relevance in handling the basic problem
of \textit{background independence}, which is inherent
to any attempt at a quantization of the theory
of gravitation.
 \smallskip

The \textit{hypoelliptic signature operator} $D$ is uniquely determined by the equation
$Q = D |D|$, where $Q$ is the operator
 \begin{equation} \nonumber
    Q = (d_V^* \, d_V - d_V \, d_V^*) \op \g_V \,  (d_H + d_H^*) \, ,
\end{equation}
acting on the Hilbert space of $L^2$-sections
\begin{equation} \nonumber
 {\FH}_{PM} = L^{2}({\wedge}^{\cdot} \Vc^{*}  \ot 
 {\wedge}^{\cdot} \Nc^{*}  , \,  \varpi_{P}) \, ; 
\end{equation}
here $d_V$ denotes the vertical exterior derivative,
$\g_V$ is the usual grading for the vertical signature operator,
$d_H $ stands for the horizontal exterior differentiation with
respect to a fixed connection on the frame bundle, and 
$\varpi_P$ is the $\Diff^+ (M)$-invariant volume form on $PM$ 
associated to the connection.
When $n \equiv 1 \, {\rm or} \, 2 \, ({\rm mod} \,  4)$, 
for the vertical component to make sense,
one has to replace $PM$ with $ PM \times S^1$ 
so that the dimension of the vertical fiber
be even. 
\medskip

The above construction allows to associate to any transversely oriented
foliation  $\Fc$ of a manifold $V$ a spectral triple encoding the geometry 
of $V/\Fc$ in the following sense.
If $M$ is a complete transversal and 
$\G$ is the corresponding holonomy pseudogroup, then
the pair $({\FH}_{PM}, D)$ described above can be completed to a 
 spectral triple $(\Ac_{\G} , {\FH}_{PM}, D)$, 
 by taking as $\Ac_{\G} $ the convolution algebra of the smooth \`etale
 groupoid associated to $\G$.  The spectral triple
$(\Ac_{\G} , {\FH}_{PM}, D)$ represents the
desired geometric structure  for  $V/\Fc$.	
\medskip

Using hypoelliptic calculus, which in particular provides 
a noncommutative residue functional $\displaystyle {\int \!\!\!\!\!\! -} $
extending the Dixmier trace,
we proved in \cite[Part I]{CMlif} that
such a spectral triple $(\Ac_{\G} , {\FH}_{PM}, D)$
fulfills the hypotheses of the
operator theoretic local index theorem of \cite[Part II]{CMlif}.
Therefore, its character-index
$ {ch}_* (D) \in H C_{*} (\Ac_{\G})$  can be expressed
in terms of residues of spectrally defined zeta-functions,
and is given by a cocycle $\{\phi_q\}$ in the
$(b,B)$ bi-complex of $\Ac_{\G}$ whose components are of the 
following form 
\begin{equation} \label{index}
\phi_q (a^0 ,\ldots ,a^q) = \sum_{\bf k} c_{q, {\bf k} } \,
{\int \!\!\!\!\!\! -}  a^0 [Q, a^1]^{(k_1)} \ldots
[Q, a^q]^{(k_q)}  \, \vert Q \vert^{-q-2\vert k\vert} \, ;
\end{equation}
we have used here the abbreviations
$T^{(k)} = \nb^k (T) \,$  and $\, \nb (T) = D^2 T - TD^2 $, 
\begin{eqnarray*}
&\ & {\bf k} \, =\, (k_1 , \ldots , k_q) \, , \qquad  \vert {\bf k}  \vert = k_1 +\ldots + k_q \, , 
\quad \text{and}   \cr \cr
 c_{q, {\bf k} }&=&
\frac{(-1)^{\vert {\bf k}  \vert}  \sqrt{2i}} {k_1 ! \ldots k_q! \, 
(k_1 +1) \ldots (k_1 + \ldots + k_q
+q)} \,  \G \left( \vert k \vert + \frac{q}{2} \right).
\end{eqnarray*}
The summation necessarily involves  finitely many nonzero terms for each
$\phi_q $, and $q$ cannot exceed $\displaystyle
\frac{n(n+1)}{2}+ 2n$. 

In practice, the actual computation of the expression (\ref{index})
is exceedingly difficult to perform. 
Even in the case of codimension $n=1$, when  
there are only two components
$\{\phi_1, \phi_3\}$, the order of magnitude of the number 
of terms one needs to handle is $10^3$.
Thus, a direct evaluation of (\ref{index}) 
 for an arbitrary codimension $n$ is impractical.
\medskip

There are two reduction steps which  help alleviate, 
to some extent, the complexity of the problem. 
First, by enlarging if necessary the
pseudogroup $\G$, one may assume that
$M$ is a flat affine manifold.  There is no loss of generality in making
this assumption as long as
 the affine structure is not required to be
preserved by $\Gamma$. Thus, one can equip $M$ with a flat
connection,  and since the horizontal component of the operator
 $Q$ is built out of the connection, its expression gets simplified
to the fullest extent possible. It is also important to note that 
$Q$ is affiliated with the universal enveloping algebra
of the group of affine motions of $\Rb^n$, in the sense that
it is of the form
 \begin{equation*} 
    Q = R (Q_{\rm alg}), \quad \text{with} \quad
   Q_{\rm alg} \in \left({\FA} (\Rb^n \rtimes {\Fg \Fl} (n, \Rb)) 
   \ot {\rm End} (E) \right)^{SO(n)}  \, ,
\end{equation*}
where $R$ is the right regular representation of $\Rb^n \rtimes GL (n, \Rb)$
and $E$ is a unitary 
$SO(n)$--module. 
\smallskip

Secondly, one can afford to work at the level of the principal
bundle $F^{+}M$, since the descent to
the quotient bundle $PM$ only
involves the simple operation
of taking $SO(n)$-invariants.
\medskip

The strategy that led to the unwinding of the formula (\ref{index}) 
essentially evolved from the following observation. The built-in 
affine invariance of the operator
$Q$, allows to reduce the noncommutative residue functional
involved in the cochains
\begin{equation}
\phi (a^0 , \ldots , a^q) = \, \int \!\!\!\!\!\! - \, a^0 [Q ,a^1]^{(k_1)}
\ldots
[Q , a^q]^{(k_q)} \, |Q|^{- (q + 2 \vert {\bf k} \vert)} \, ,  
\end{equation}
to a genuine integration, and thus replace them by cochains
of the form
\begin{equation} \label{tc}
\psi (a^0 , \ldots , a^q) = \tau_{\G} (a^0 \, h^1 (a^1) 
\ldots h^q (a^q)) \, ;
\end{equation}
here $\tau_{\G}$ is the canonical trace on $\Ac_{\G}$ and 
$  h^1 , \ldots ,  h^q $ are `transverse' differential operators 
acting on the algebra $\Ac_{\G}$. Under closer scrutiny, which will be
discussed in great detail in the next section, these transverse differential operators 
turn out to arise from the action
of a canonical Hopf algebra $\Hc_n$, depending only on the
codimension $n$. Furthermore, 
the cochains (\ref{tc}) will be recognized to belong to
the range of a certain cohomological characteristic map. 
\bigskip

\section{The Hopf algebra $\Hc_n$ and its standard action}

Let  $F{\Rb}^n$
be the frame bundle on ${\Rb}^n$, identified to
${\Rb}^n \times GL(n, {\Rb})$ in the usual way:  the 1-jet at $0 \in {\Rb}^n$ of the map
$\, \phi : {\Rb}^n \ra {\Rb}^n$,
$$
\phi(t) = x \, + \, {\bf y}  t  \, , \qquad x, \, t \in \Rb^n \, , \quad {\bf y}  \in GL(n, {\Rb})
$$
is identified to the pair $(x, {\bf y}) \in GL(n, {\Rb})$. We endow it with the trivial
connection, given by the matrix-valued $1$-form $\, \om = ( \om^i_j ) \,$ where,
with the usual summation convention,
\begin{equation} \label{trivcon}
 \om^i_j \, :=  \, ({\bf y}^{-1})^i_{\mu} \, d{\bf y}^{\mu}_j \, = \, ({\bf y}^{-1} \, d{\bf y})^i_j
 \, , \qquad i, j =1, \ldots , n \, .
 \end{equation}
The corresponding basic horizontal vector fields on $F{\Rb}^n$ are
\begin{equation}
 X_k \, = \, y_k^{\mu} \, \part_{\mu} \, , \quad   k =1, \ldots , n \, ,
 \quad \text{where} \quad \part_{\mu} = \frac{\part} {\part \, x^{\mu}}  \, .
 \end{equation}
We denote by $\, \t = ( \t^k ) \,$ the canonical form of the frame bundle
\begin{equation}  \label{can}
 \t^k \, :=  \, ({\bf y}^{-1})^k_{\mu} \, dx^{\mu} \, = \, ({\bf y}^{-1} \, dx)^k
 \, , \qquad k =1, \ldots , n \,
 \end{equation}
 and then let
\begin{equation}
 Y_i^j \, = \,  y_i^{\mu} \, \part_{\mu}^j \, , \quad  i, j =1, \ldots , n \, ,
 \quad \text{where}
\quad \part_{\mu}^j  := \frac{\part}{\part \, y_j^{\mu}}  \, ,
 \end{equation}
be the fundamental vertical vector fields associated
to the standard basis of ${\Fg \Fl} (n, \Rb )$ and
generating the canonical right action
of $GL(n, {\Rb})$ on $F{\Rb}^n$.
At each point of $F{\Rb}^n$,
$\{ X_k , Y_i^j \} $, resp.  $\{ \t^k , \om^i_j \}$, form bases of the tangent, resp.
cotagent space, dual to each other:
\begin{eqnarray} \label{dual}
 \langle \om^i_j , Y_k^\ell \rangle \, &=& \, \d^i_k \d^\ell_j \, , \quad \langle \om^i_j , X_k \rangle \, = \, 0 \, , \\ \nonumber
 \langle \t^i, Y_k^\ell \rangle \, &=& \, 0 \, , \qquad \langle \t^i, X_j \rangle \, = \, \d^i_j \, .
 \end{eqnarray}
\medskip

The group of diffeomorphisms
 $ \Gc_n := \Diff {\Rb}^n$
acts on $F{\Rb}^n$, by the natural lift of the tautological action to the frame level:
\begin{equation}
\wt{\vp} (x, {\bf y}) := \left( \vp (x), {\vp}^{\prime} (x) \cdot {\bf y} \right) \, , \quad \text{where}
\quad {\vp}^{\prime} (x)^{i}_{j} = \part_{j} \, {\bf \vp}^{i}  (x) \, .
\end{equation}
Viewing here $\Gc_n$ as a discrete group,
we  form the crossed product algebra
$$
   {\Ac}_n \, : = \, C_c^{\ify} (F{\Rb}^n ) \semi \Gc_n  \, .
$$
As a vector space, it is spanned
by monomials of the form $\,f \, U_{\vp}^* \,$, where
$\, f \in C_c^{\ify} (F{\Rb}^n) \, $ and  $\, U_{\vp}^* \,$
stands for $\, \wt{\vp}^{-1} $, while the product is given by the
multiplication rule
\begin{equation}
f_1 \, U_{\vp_1}^* \cdot f_2 \, U_{\vp_2}^* =
f_1 (f_2 \circ \wt{\vp}_1) \,
U_{\vp_2 \vp_1}^* \, .
\end{equation}
Alternatively, $\Ac_n$ can be regarded as the subalgebra of the endomorphism algebra
$\Lc \left(C_c^{\ify} (F{\Rb}^n ) \right)$ of the vector space $  \, C_c^{\ify} (F{\Rb}^n ) \, $,
generated by the multiplication and the translation operators
\begin{eqnarray}
M_f (\xi)   &=& \,  f \, \xi \, , \quad f \in C_c^{\ify} (F{\Rb}^n) \, , \, \xi \in C_c^{\ify} (F{\Rb}^n ) \\
U_{\vp}^* (\xi)  &=& \,  \xi \circ \wt{\vp}  \, , \qquad \vp \in \Gc_n \, , \,
\xi \in C_c^{\ify} (F{\Rb}^n ) \, .
\end{eqnarray}

Since the right action of $GL(n,\Rb)$ on  $F{\Rb}^n$ commutes with the action of
$\Gc_n$, at the Lie algebra level one has
\begin{equation} \label{Ydisp}
U_{\vp} \, Y_i^j  \, U_{\vp}^* \, = \, Y_i^j  \, , \qquad  \vp  \in \Gc_n \, .
\end{equation}
This allows to promote the vertical vector fields to derivations of $\Ac_n$.
Indeed, setting
\begin{equation}
Y_i^j (f \, U_{\vp}^*) \, = \, Y_i^j ( f) \, U_{\vp}^*  \, , \quad
f \, U_{\vp}^* \in {\Ac}_n \, ,
\end{equation}
the extended operators satisfy the derivation rule
\begin{equation} \label{Yrule}
Y_i^j (a \, b) \, = \, Y_i^j (a) \, b \, + \, a \,Y_i^j (b) \,   , \quad
a, b \in {\Ac}_n \, ,
\end{equation}
\medskip

We shall also prolong the horizontal vector fields to
linear transformations $X_k  \in \Lc \, ({\Ac_n}) $, in a
similar fashion:
\begin{equation}
X_k (f \, U_{\vp}^* ) = X_k (f) \, U_{\vp}^* \, , \quad
f \, U_{\vp}^* \in {\Ac}_n \, .
\end{equation}
The resulting operators are no longer $\Gc_n$-invariant.
Instead of (\ref{Ydisp}), they satisfy
\begin{equation} \label{Xdisp}
 U_{\vp} \, X_k \, U_{\vp}^* \,  =  \,  X_{k} \,  - \, \g_{jk}^i (\vp^{-1}) \, Y_i^j  \, ,
\end{equation}
where $ \, \vp \mapsto \g_{jk}^i (\vp) \,$ is a group $1$-cocycle on $\Gc_n$
with values in $C^{\ify} (F{\Rb}^n ) $; specifically,
\begin{equation} \label{gijk}
\g_{jk}^i (\vp) (x, {\bf y}) \, =\, \left( {\bf y}^{-1} \cdot {\vp}^{\prime} (x)^{-1} \cdot \part_{\mu} {\vp}^{\prime} (x)
\cdot {\bf y}\right)^i_j \, {\bf y}^{\mu}_k \, .
\end{equation}
The above expression comes out readily from the pull-back formula for the connection,
\begin{equation} \label{omdisp}
\wt{\vp}^* (\om^i_j ) \, = \, \om^i_j  \, + \, \g_{jk}^i (\vp)  \,\t ^k  \, ;
\end{equation}
indeed, if $\, ( \wt{x}, \wt{\bf y} ) := \wt{\vp} (x, {\bf y}) = \left( \vp (x), {\vp}^{\prime} (x) \cdot {\bf y} \right)$, then
\begin{eqnarray*}
  {\wt{\bf y} }^{-1} \, d \wt{\bf y} \, &=& \, {\bf y}^{-1} \, {\vp}^{\prime} (x)^{-1} \left( d{\vp}^{\prime} (x) \, {\bf y} \, + \,
  {\vp}^{\prime} (x) \, d{\bf y}\right) \\ \nonumber
   &=& \, {\bf y}^{-1} \, d{\bf y} \, + \,
  \left({\bf y}^{-1} \, {\vp}^{\prime} (x)^{-1} \, \part_{\mu} {\vp}^{\prime} (x) \, {\bf y}\right) \, dx^{\mu} \\ \nonumber
  &=& \, {\bf y}^{-1} \, d{\bf y} \, + \,
 \left({\bf y}^{-1} \, {\vp}^{\prime} (x)^{-1} \, \part_{\mu} {\vp}^{\prime} (x) \, {\bf y}\right) \,{\bf y}^{\mu}_k \, \t ^k  \, .
\end{eqnarray*}
In view of the $\Gc_n$-invariance of $\t$, (\ref{omdisp}) makes the cocycle property of
$\g$ obvious. To obtain (\ref{Xdisp}), one just has to use that   $\{ \t^k , \, (\wt{\vp}^{-1})^* (\om_i^j ) \}$
is the dual basis to $\{U_{\vp} \, X_k \, U_{\vp}^*  , \, Y_i^j \} $, cf. (\ref{dual}).
\medskip

As a consequence of (\ref{Xdisp}), the operators $X_k  \in \Lc \, ({\Ac_n}) $
are no longer derivations of $\Ac_n$, but satisfy
instead a non-symmetric Leibniz rule:
\begin{equation} \label{Xrule}
X_k(a \, b) \, = \, X_k (a) \, b \, + \, a \,X_k (b) \, + \, \d_{jk}^i (a) \, Y_{i}^{j} (b)
 \,   , \quad
a, b \in {\Ac}_n \, ,
\end{equation}
where the linear operators $\, \d_{jk}^i   \in \Lc \, ({\Ac_n}) $ are defined by
\begin{equation} \label{d}
\d_{jk}^i (f \, U_{\vp}^*) \, =\, \g_{jk}^i  (\vp) \, f \, U_{\vp}^*  \, .
\end{equation}
Indeed, on taking $\, a=f_1 \, U_{\vp_1}^*$, $\, b = f_2 \, U_{\vp_2}^*$,  one has
\begin{eqnarray*}
X_k (a \cdot b) &=&  X_k (f_1 \, U_{\vp_1}^* \cdot f_2 \, U_{\vp_2}^*) \,
 = \, X_k (f_1 \cdot U_{\vp_1}^* \, f_2 \, U_{\vp_1}) \, U_{\vp_2 \vp_1}^* \\ \nonumber
&=& X_k (f_1) \, U_{\vp_1}^* \cdot f_2 \, U_{\vp_2}^* \,
+ \, f_1 \, U_{\vp_1}^* \cdot X_k (f_2 \, U_{\vp_2}^*)  \\ \nonumber
&+& f_1 \, U_{\vp_1}^* \cdot
(U_{\vp_1} \, X_k \,  \, U_{\vp_1}^* \,  \, - \, X_k) ( f_2 \, U_{\vp_2}^*) \, ,
\end{eqnarray*}
which together with (\ref{Xdisp}) and the cocycle property of $\g_{jk}^i$ imply
(\ref{Xrule}).

The same cocycle property
shows that the operators $\, \d_{jk}^i $ are derivations:
\begin{equation} \label{drule}
\d_{jk}^i (a \, b) \, = \, \d_{jk}^i  (a) \, b \, + \, a \, \d_{jk}^i  (b) \,   , \quad
a, b \in {\Ac}_n \, ,
\end{equation}

The operators $\, \{X_k , \, Y^i_j \} \,$
satisfy the commutation relations of the group of affine transformations of ${\Rb}^n$:
\begin{eqnarray} \label{aff}
[Y_i^j , Y_k^{\ell}] &=& \d_k^j Y_i^{\ell} - \d_i^{\ell} Y_k^j \, , \\ \nonumber
[Y_i^j , X_k] &=& \d_k^j X_i \, , \qquad \qquad
[X_k , X_{\ell}]\, =\, 0 \, .  \nonumber
\end{eqnarray}
The successive commutators of the operators $\, \d_{jk}^i $'s
with the $X_{\ell}$'s yield  new generations of
\begin{equation} \label{highd}
\d_{jk \vert \ell_1 \ldots \ell_r}^i \, := \,
[X_{\ell_r} , \ldots [X_{\ell_1} , \d_{jk}^i] \ldots ] \, ,
\end{equation}
which involve multiplication by higher order jets of diffeomorphisms
\begin{eqnarray} \label{d'}
\d_{jk \vert \ell_1 \ldots \ell_r}^i \, ( f \, U_{\vp}^*) \, &=& \,
 \g_{jk \vert \ell_1 \ldots \ell_r}^i \, f \, U_{\vp}^* \, , \qquad \text{where}  \\
  \g_{jk \vert \ell_1 \ldots \ell_r}^i \, &=&
\, X_{\ell_r} \cdots X_{\ell_1} (\g_{jk}^i) \, . \nonumber
\end{eqnarray}
Evidently, they commute among themselves:
\begin{equation} \label{abel}
[\d_{jk \vert \ell_1 \ldots \ell_r}^i , \, \d_{j'k' \vert {\ell}'_1 \ldots  {\ell}_r}^{i'}] \, = \, 0 \, .
\end{equation}
The operators $\, \d_{jk \vert \ell_1 \ldots \ell_r}^i \,$
are not all distinct; the order of the first two lower
indices or of the last $r$ indices is immaterial. Indeed,
performing the matrix multiplication in (\ref{gijk}) gives the expression
\begin{equation} \label{sgijk}
\g_{jk}^i (\vp) (x, {\bf y}) \, =\, ( {\bf y}^{-1} )^i_\lb \, ({\vp}^{\prime} (x)^{-1})^\lb_\rho
 \, \part_{\mu} \part_\nu {\vp}^{\rho} (x) \, {\bf y}^\nu_j \, {\bf y}^{\mu}_k \, ,
\end{equation}
which is clearly symmetric in the indices $j$ and $k$. The symmetry in
the last $r$ indices follows from the definition (\ref{highd}) and the fact that,
the connection being flat, the horizontal vector fields commute. It can also be
directly seen from the explicit formula
\begin{eqnarray} \label{highg}
&\g_{jk \vert \ell_1 \ldots \ell_r}^i (\vp) (x, {\bf y}) \, = \\ \nonumber
 &= \, ({\bf y}^{-1})^i_\lb \part_{\b_r} \ldots \part_{\b_1}
\left(({\vp}^{\prime} (x)^{-1})^\lb_\rho
 \part_{\mu} \part_\nu {\vp}^{\rho} (x) \right) {\bf y}^\nu_j  {\bf y}^{\mu}_k  {\bf y}^{\b_1}_{\ell_1}
 \ldots {\bf y}^{\b_r}_{\ell_r}\, .
\end{eqnarray}
\bigskip

The commutators between the $Y_\nu^\lb$'s and $\d_{j k}^i$'s can be obtained
from the explicit expression (\ref{gijk}) of the cocycle $\g$ ,
by computing its derivatives
in the direction of the vertical vector fields. Denoting by $E_\nu^\lb$ the $n \times n$
matrix whose entry at the $\lb$-th row and the $\mu$-th column is $1$ and all others
are $0$, one has:
\medskip

$\, Y_\nu^\lb (\g_{jk}^i (\vp) ) (x, {\bf y}) \, = $
\begin{eqnarray*}
 &=& \frac{d}{dt} \vert_{t=0} \left(
\left(\exp(-t E_\nu^\lb) \, {\bf y}^{-1} \, {\vp}^{\prime} (x)^{-1} \, \part_{\mu} {\vp}^{\prime} (x)
\, {\bf y} \exp(t E_\nu^\lb) \right)^i_j \, ({\bf y} \exp(t E_\nu^\lb))^{\mu}_k \right) \, \\ \nonumber
 &=&  - [E_\nu^\lb , \, {\bf y}^{-1} \, {\vp}^{\prime} (x)^{-1} \, \part_{\mu} {\vp}^{\prime} (x)
\, {\bf y} ]^i_j \, {\bf y}^{\mu}_k \, + \,
\left({\bf y}^{-1} \, {\vp}^{\prime} (x)^{-1} \, \part_{\mu} {\vp}^{\prime} (x)
\, {\bf y} \right)^i_j \, ({\bf y} E_\nu^\lb)^{\mu}_k \\ \nonumber
&=& \left( {\bf y}^{-1} \, {\vp}^{\prime} (x)^{-1} \, \part_{\mu} {\vp}^{\prime} (x)
\, {\bf y}\right)^i_{\nu} \, \d_j^{\lb}\, {\bf y}^{\mu}_k \, - \, \d^i_\nu \,
 \left( {\bf y}^{-1} \, {\vp}^{\prime} (x)^{-1} \, \part_{\mu} {\vp}^{\prime} (x)
\, {\bf y} \right)^\lb_j \, {\bf y}^{\mu}_k \\ \nonumber
&+& \,
\left({\bf y}^{-1} \, {\vp}^{\prime} (x)^{-1} \, \part_{\mu} {\vp}^{\prime} (x)
{\bf y} \right)^i_j \, {\bf y}^{\mu}_{\nu} \, \d_k^{\lb} \, .
\end{eqnarray*}
It follows that
\begin{equation*} \label{Yd}
 [Y_\nu^\lb , \d_{j k}^i ] \, =\,  \d_j^\lb \ \d_{\nu k}^i  \, + \, \d_k^\lb \ \d_{j \nu}^i
 \, - \, \d_\nu^i \ \d_{j k}^\lb  \, .
\end{equation*}
More generally, using (\ref{aff}) and the very definition (\ref{highd}), one
obtains from this, by induction,
\begin{equation}  \label{Yhighd}
 [Y_\nu^\lb , \d_{j_1 j_2 \vert j_3 \ldots j_r}^i ] \, =\, \sum_{s=1}^r \, \d^\lb_{j_s} \,
 \d^i_{j_1  j_2 \vert j_3  \ldots j_{s-1} \nu j_{s+1} \ldots j_r}
 \, - \, \d_\nu^i \, \d_{j_1 j_2 \vert j_3  \ldots j_r}^\lb  \, .
\end{equation}
\bigskip

By a transient abuse of notation,
we now regard $X_k $, $Y^i_j$ and  $\d_{jk \vert \ell_1 \ldots \ell_r}^i $
as abstract symbols, preserving the convention that in the designation of the latter
the order of the first two lower
indices or of the last $r$ indices is unimportant, and make
the following definition.

\begin{definition} \label{def}
Let  $\Hc_n$ be the universal enveloping algebra of the Lie
algebra  $ {\Fh}_n$  with basis
$$ \{X_\lb , \, Y^\mu_\nu , \, \d_{jk \vert \ell_1 \ldots \ell_r}^i
\vert 1 \leq \lb, \mu, \nu , i \leq n, \, 1\leq j \leq k \leq n, \, 1 \leq \ell_1\leq \ldots \leq \ell_r \leq n   \}
$$
and the following presentation:
\begin{eqnarray} \label{Xcom}
[X_k , X_{\ell}] &=& 0 , \\ \label{Ycom}
[Y_i^j , Y_k^{\ell}] &=& \d_k^j Y_i^{\ell} - \d_i^{\ell} Y_k^j  , \\ \label{YXcom}
[Y_i^j , X_k] &=& \d_k^j X_i  , \\\label{Xdcom}
[X_{\ell_r} , \d_{jk \vert \ell_1 \ldots \ell_{r-1}}^i ] &=&  \d_{jk \vert \ell_1 \ldots \ell_r}^i  
 ,\\ \label{Ydcom}
[Y_\nu^\lb , \d_{j_1 j_2 \vert j_3 \ldots j_r}^i ] &=& \sum_{s=1}^r \, \d^\lb_{j_s} \,
 \d^i_{j_1  j_2 \vert j_3  \ldots j_{s-1} \nu j_{s+1} \ldots j_r}
  -  \d_\nu^i \, \d_{j_1 j_2 \vert j_3  \ldots j_r}^\lb , \\ \label{ddcom}
[\d_{jk \vert \ell_1 \ldots \ell_r}^i , \, \d_{j'k' \vert \ell'_1 \ldots \ell'_r}^{i'}]  &=& \, 0 \, .  
\end{eqnarray}
\end{definition}
\bigskip

We shall endow  $\, \Hc_n = \FA ({\Fh}_n)\,$ with a canonical Hopf structure, 
which is non-cocommutative and therefore different from the standard structure of
a universal enveloping algebra.
\medskip

\begin{proposition} \label{hopf} 
$1^0$. The formulae
\begin{eqnarray}  \label{DX}
\D  X_k &=& \,  X_k \ot 1 + 1\ot X_k + \d^i_{jk} \ot Y_i^j \, , \\   \label{DY}
\D Y_i^j &=& \, Y_i^j \ot 1 + 1 \ot Y_i^j \, , \\    \label{Dd}
\D \d^i_{jk} &=& \d^i_{jk} \ot 1 + 1 \ot \d^i_{jk} \, ,    
\end{eqnarray}
uniquely determine a 
coproduct  $\, \D : \Hc_n \ra \Hc_n \ot \Hc_n\,$,
which makes $\Hc_n $ a bialgebra with respect to the product 
$\, m : \Hc_n \ot \Hc_n \ra \Hc_n \,$
and the counit  $ \ve: \Hc_n \ra {\Cb}$ inherited from $ \FA ({\Fh}_n) $.
\smallskip

$\qquad 2^0$. The formulae
\begin{eqnarray}  \label{SX}
S (X_k) &=& \, - X_k  + \d^i_{jk} \, Y_i^j \, , \\  \label{SY}
S(Y_i^j) &=& \, - Y_i^j \, , \\   \label{Sd}
S(\d^i_{jk}) &=& - \d^i_{jk}  \, ,    
\end{eqnarray}
uniquely determine an anti-homomorphism  $\, S: \Hc_n \ra \Hc_n \,$,
which provides the antipode that turns $\Hc_n $ into a Hopf algebra. 
\end{proposition}
\medskip

\begin{proof}  
$1^0$. Once its existence is established, the uniqueness of the coproduct
$\D$ satisfying (\ref{DX})--(\ref{Dd}) is obvious. Indeed, these formulae
prescribe the values of
the algebra homomorphism $\D: \FA ({\Fh}_n) \ra  \FA ({\Fh}_n) \ot \FA ({\Fh}_n)$ 
on a set of generators. (Note however that, in view of
(\ref{DX}), they \textit{do not} define a
 Lie algebra homomorphism ${\Fh}_n \ra {\Fh}_n \ot  {\Fh}_n $).
 
 To prove the existence of the coproduct, one checks that the presentation
 (\ref{Xcom})--(\ref{ddcom}) is preserved by $\D$; this ensures that  $\D$
extends to an \textit{algebra homorphism}  
$\FA ({\Fh}_n) \ra  \FA ({\Fh}_n) \ot \FA ({\Fh}_n)$.
\medskip

One then has to verify \textit{coassociativity} and \textit{counitality}. We skip
the straightforward details. An alternate argument
will emerge later (see Remark \ref{copro}).
 \medskip
 
$2^0$.  Similarly, one shows that $S$ defines an an anti-homomorphism  
 $\, \Hc_n \ra \Hc_n \,$ by checking that the  presentation
 (\ref{Xcom})--(\ref{ddcom}) is anti-preserved. Since 
 the \textit{antipode} axioms  are multiplicative, it 
 suffices to verify them on a set of generators. This is precisely how
 the formulae (\ref{SX})--(\ref{Sd}) were obtained, and is easy to
 verify.
 \end{proof}
\bigskip

The abuse of notation made in the definition \ref{def}
will be rendered completely innocuous by the next result.
\medskip

\begin{proposition} \label{free}
$1^0$. The subalgebra of $\, \Lc \, ({\Ac_n}) $ generated by the linear operators
$\, \{X_k , \, Y^i_j , \, \d_{jk}^i \vert \, i, j, k =1, \ldots , n \} \,$
is isomorphic to the algebra $ \Hc_n$ .
\smallskip

$\qquad 2^0$. The ensuing action of $\Hc_n $ turns $ \Ac_n $ into a left
$\Hc_n $-module algebra.
\end{proposition}
\medskip

\begin{proof}
$1^0$. The notation needed to specify a Poincar\'e-Birkhoff-Witt basis for
$ \FA ({\Fh}_n) $ involves two kinds of multi-indices.
The first kind are of the form
$$
 I \, = \left\{ i_1 \leq \ldots \leq i_p \, ; \left( ^{j_1}_{k_1} \right) \leq \ldots \leq
 \left( ^{j_q}_{k_q} \right) \right\} \, ,
$$
while the second kind are of the form
$\, \displaystyle K \, = \, \{ \k_1 \leq \ldots \leq \k_r \} $, where
\begin{equation*}
\k_s \, = \, \left( \begin{matrix}
i_s &\qquad \qquad\qquad \hfill \cr
 j_s & k_s \,\; \vert\,
\ell_1^s \leq \ldots \leq \ell_{p_s}^s \hfill \cr\end{matrix}
\right)  \, , \qquad s=1, \ldots , r \, ; \nonumber
\end{equation*}
in both cases the inner multi-indices are ordered lexicographically.
\medskip 

The PBW basis
of $\FA ({\Fh}_n) $ will consist of elements of the form $\, \d_K \, Z_I \,$,
ordered lexicographically, where
$$
Z_I = X_{i_1} \ldots X_{i_p} \, Y_{k_1}^{j_1} \ldots Y_{k_q}^{j_q} 
\quad \hbox{and} \quad
\d_K = \d_{j_1 \, k_1 \vert  \ell_1^1 \ldots \ell_{p_1}^1}^{i_1} \ldots \d_{j_r
\, k_r \vert  \ell_1^r \ldots \ell_{p_r}^r}^{i_r} \, .
$$
\medskip

We need to prove that if $\, c_{I, \k} \in \Cb$ are such that
\begin{equation} \label{0}
 \sum_{I, K} \, c_{I, K} \,  \d_K \, Z_I   \, (a) \, =  \, 0 \, ,
\quad \fl \,  a \in \Ac_n \, , 
\end{equation}
then $\, c_{I, K} = 0$, for any $(I, K)$.

To this end, we evaluate (\ref{0}) on all monomials
$\, a \, = \, f \, U_{\vp}^* $ at the point 
$$ u_0 \, = \, (x=0, {\bf y} = {\bf I}) \in F {\Rb}^n = {\Rb}^n \times GL(n,  {\Rb}^n) \, .
$$ 
In particular, for any fixed but arbitrary $\vp \in \Gc_n$, one obtains
\begin{equation}  \label{00}    
\sum_{I} \left(\sum_{K} \, c_{I, K}  \, \g_{\k} (\vp) ( u_0)\right)
   (Z_I \, f) ( u_0) = 0 \,  , \quad \fl \, f \in C_c^{\ify} (F{\Rb}^n ) \, .
\end{equation} 
Since the $Z_I$'s form a PBW basis of $\FA ({\Rb}^n \semi {\Fg \Fl} (n, \Rb)) $,
which can be viewed as
the algebra of left-invariant differential operators on $F{\Rb}^n $,
the validity of (\ref{00}) for any $f \in C_c^{\ify} (F{\Rb}^n )$
implies the vanishing for each $I$ of the corresponding coefficient. 
One therefore obtains, for any fixed $I$,
\begin{equation}  \label{01}    
\sum_{K} \, c_{I, K}  \, \g_{K}  (\vp) ( u_0) \,
 \, = 0 \,  , \quad \fl \, \vp \in \Gc_n \, .
\end{equation} 
To prove the vanishing of all the coefficients, 
we shall use induction on the \textit{height} of 
$\, \displaystyle K \, = \, \{ \k_1 \leq \ldots \leq \k_r \} $;
the latter is defined by counting the total
number of horizontal derivatives of its largest
components:
$$ 
\vert K \vert \, = \, \ell^r_1 + \cdots + \ell^r_{p_r} \, .
$$
\medskip

We start with the case of height $0$, when the identity (\ref{01}) reads
\begin{equation*}    
\sum_{K} \, c_{I, K}  \, \g_{j_1 k_1}^{i_1} (\vp)( u_0)  \cdots 
 \g_{j_r k_r}^{i_r} (\vp) ( u_0 ) \,
 \, = 0 \,  , \quad \fl \, \vp \in \Gc_n \, .
\end{equation*} 
 Choosing $\vp$  in the subgroup $\Gc_n^{(2)} (0) \subset \Gc_n$
consisting of the diffeomorphisms  whose 2-jet at $0$ is of the form
\begin{equation*}    
J_0^2 ( \vp)^i (x) \, = \, x^i \, + \,
\frac{1}{2}  \sum_{j, k =1}^n \xi^{i}_{j k} x^j x^k, 
\quad  \xi \in \Rb^{n^{3}}, \quad \xi^{i}_{j k} =   \xi^{i}_{k j} \, ,
\end{equation*} 
and using (\ref{sgijk}), one obtains:
\begin{equation*}
\sum_{K} \, c_{I, K}  \, \xi_{j_1 k_1}^{i_1}  \cdots 
 \xi_{j_r k_r}^{i_r} \, = 0 , \quad  
\xi^{i}_{j k}  \in \Rb^{n^{3}}, \quad  \xi^{i}_{j k} =   \xi^{i}_{k j}  \, .
\end{equation*} 
It follows that all coefficients $\, c_{I, K} \, = \, 0$.
\medskip

Let now $N \in \Nb$ be the largest height of  occurring in (\ref{01}).
By varying $\vp $ in the 
subgroup $\Gc_n^{(N+2)} (0) \subset \Gc_n$ of
all diffeomorphisms whose $(N+2)$-jet at $0$ has the form
\begin{eqnarray*}    
J_0^{N+2}  (\vp)^i (x) &=  \, x^i +
\frac{1}{(N+2)!}  \sum_{j, k, \a_1, \ldots , \a_{N+2} } \xi^i_{j k \a_1 \ldots \a_{N}}  
x^j x^k x^{\a_1} \cdots x^{\a_{N}} ,  \\ \nonumber
  \xi^i_{j k \a_1 \ldots \a_{N}} & \in \Cb^{n^{N+3}}  ,  \, \,
 \xi^i_{j k \a_1 \ldots \a_{N}} = \xi^i_{k j \a_{\s(1)} \ldots \a_{\s(N)}} ,
 \, \, \fl \, \, \text{permutation} \, \, \s \, , \nonumber
\end{eqnarray*} 
and using (\ref{highg}) instead of (\ref{sgijk}), one derives as above  
 the vanishing
of all coefficients $c_{I, \k}  $ with $\vert \k \vert = N$. This lowers the
height  in (\ref{01})  and thus completes the induction.
\medskip

$2^0$. Due to its multiplicative nature, it suffices to check the 
compatibility property
\begin{equation} \label{HA}
 h (a  b) \ = \ \sum_{(h)} \, h_{(1)} (a) \, h_{(2)} (b) \, , 
\qquad  h \in \Hc_n \, , \quad a, b \in \Ac_n 
\end{equation}
only on generators. This is precisely what the formulae (\ref{Yrule}), (\ref{Xrule})
and (\ref{drule}) verify.
  \end{proof}
\bigskip

As a matter of fact, the action of the algebra $\Hc_n$ on  $\Ac_n$ is not only
faithful but even \textit{multi-faithful}, in the sense made clear by the 
result that follows. In order to state it, we associate to an element
$h^1 \ot \ldots \ot  h^p \in \Hc_n^{ \ot ^p} $ a
\textit{multi-differential} operator acting on $\Ac_n$, as follows:
\begin{eqnarray} \label{T}
T(h^1 \ot \ldots \ot  h^p) \, (a^1\ot \ldots \ot  a^p) 
\, = \, h^1(a^1) \cdots h^p (a^p) \, , \\
\text{where}  \quad  h^1, \ldots , h^p \in \Hc_n \, \quad \text{and} \quad
a^1, \ldots , a^p \in\Ac_n \, .\nonumber
\end{eqnarray}
\smallskip

\begin{proposition} \label{multid}
The linearization $ T\, : \, \Hc_n^{ \ot ^p}  \longra \, \Lc (\Ac_n^{ \ot^ p} , \Ac_n) \, $
of the above assignment  is injective for each $p \in \Nb$.
\end{proposition}
\smallskip

\begin{proof}  For $p=1$, $T$ gives the standard action of 
$\Hc_n$ on $\Ac_n$, which was just shown to be faithful.
To prove that $ \Ker T \, = \, 0 \, $ for an arbitrary $p \in \Nb$,
assume that
\begin{equation*} 
H \, =\, \sum_{\rho} \, h^{1}_{\rho} \ot \cdots \ot h^{p}_{\rho} \, \in 
 \Ker T  \, .
\end{equation*}
After fixing a Poincar\'e-Birkhoff-Witt basis  as above,
we may uniquely express each $h^{j}_{\rho}$ in the form
$$
h^j_{\rho} = \, \sum_{I_{j}, K_{j}} \ C_{\rho , \, I_{j}, K_{j}}
 \,  \d_{K_{j}} \, Z_{I_{j} } \, , \quad
\text{with} \quad C_{\rho , \, I_{j}, K_{j}} \in \Cb \, .
$$
Evaluating $T(H)$ on elementary tensors of the form
$ \, f_1 U_{\vp_1}^* \ot \cdots  \ot f_p U_{\vp_p}^*  \, $, 
one obtains
\begin{equation*}
\sum_{\rho, I , K} \, C_{\rho, \, I_{1}, K_{1}} \cdots C_{\rho, \, I_{p}, K_{p}} 
 \,  \d_{K_{1}} \left(Z_{I_{1} } (f_1) U_{\vp_1}^* \right)  \cdots \,
\d_{K_{p}} \left( Z_{I_{p} } (f_p) U_{\vp_p}^* \right) \, = \, 0 \, .
\end{equation*}
Evaluating further at a point $\, u_1 = (x_1, {\bf y}_1) \in F {\Rb}^n$,
and denoting
$$
u_2 = \wt{\vp}_1 (u_1) \, ,  \ldots , \, u_{p} = 
\wt{\vp}_{p-1} (u_{p-1}) \, ,
$$
the above identity gives
\begin{eqnarray*}
\sum_{\rho, I , K} \, C_{\rho, \, I_{1}, K_{1}} \cdots C_{\rho, \, I_{p}, K_{p}} 
&\cdot&  \g_{K_{1}} ({\vp}_1) (u_1)  \cdots \g_{K_{p}}  ({\vp}_p) (u_p) \\ \nonumber
&\cdot&  Z_{I_{1} } (f_1) (u_1) \cdots   Z_{I_{p} } (f_p) (u_p)\, = \, 0 \, . \nonumber
\end{eqnarray*}

Let us fix points $\, u_1, \, \ldots , \, u_p \in  F {\Rb}^n$ and then
diffeomorphisms
$\psi_0, \, \psi_1, \, \ldots , \, \psi_p$, such that
$$  u_2 = \wt{\psi}_1 (u_1) \, ,  \ldots , \, u_{p} = 
\wt{\psi}_{p-1} (u_{p-1}) \, .
$$
Following a line of reasoning similar to that of the preceeding proof, 
and iterated with respect to the points $\, u_1, \ldots , u_p \,$, 
we can infer that for  each $p$-tuple of indices of the first kind
$ \, ( I_1 , \ldots , I_p) \, $ one has
\begin{equation*}
\sum_{\rho, K} \, C_{\rho, \, I_{1}, K_{1}} \cdots C_{\rho, \, I_{p}, K_{p}} 
\cdot \g_{K_{1}} ({\vp}_1) (u_1)  \cdots \g_{K_{p}}  ({\vp}_p) (u_p) 
 \, = \, 0 \, .
\end{equation*}
Similarly, making repeated use of diffeomorphisms of the form 
$$
\psi_k  \circ \vp \quad \text{with} \quad \vp \in  \Gc_n^{(N)} (u_k) \, , \quad 
k=\, 1, \ldots, p \, ,
$$
for sufficiently many values of $N$,
we can eventually conclude that for any $ \, ( K_1 , \ldots , K_p) \,$
\begin{equation*} \label{nul*}
\sum_{\rho} \, C_{\rho, \, I_{1}, K_{1}} \cdots C_{\rho, \, I_{p}, K_{p}} 
 \, = \, 0 \, .
\end{equation*}
This proves that $\, H \, = \, 0$.
\end{proof}
\bigskip

\begin{remark}  \label{copro}
{\rm The coproduct $\D: \Hc_n \ra \Hc_n \ot \Hc_n$ and its fundamental properties
are completely determined by the 
action of $\, \Hc_n$ on its} standard module algebra $\Ac_n$.
\end{remark}
Indeed, the compatibility rule (\ref{HA}) 
can be rewritten as
\begin{equation} \label{SA}
 T(\D h) (a\ot b)\, = \,  h (a  b)  \, , \quad   \fl \,
 h \in \Hc_n \, , \,  a, b \in \Ac_n \, .
 \end{equation}
By Proposition \ref{multid} this completely determines $\D$.
Furthermore,
the coassociativity of $\D$ becomes a consequence of the associativity
of $\Ac_n$, because after applying $T$ it amounts to the identity
\begin{equation*} 
h ( (a b) c) \, = \,  h ( a (b c)) \, , \quad \fl \, h \in \Hc_n \, , \,  a, b \in \Ac_n \, .
 \end{equation*}
Similarly, the counitality is a byproduct of the unitality of  $\Ac$; transported
via $T$ it becomes tantamount to
\begin{equation*} 
h ( (a \, 1) \, = \,  h ( 1 \, a) \, = \,h(a) \, , 
\quad   \fl \, h \in \Hc_n \, , \,  a \in \Ac_n \, .
 \end{equation*} 
\bigskip

\section{Invariant trace and characteristic cochains}

An important feature of the standard module algebra 
is that it carries an \textit{invariant trace}, 
uniquely determined up to a
scaling factor. It is defined as the linear functional $\tau:\Ac_n \ra \Cb$,
\begin{equation} \label{tr} 
\tau \, (f \, U_{\vp}^* )  \, = \, 
\left\{ \begin{matrix} 
  \int_{F{\Rb}^n} \, f  \, \varpi \, , \quad \text{if}
\quad \vp = Id \, , \cr\cr
\quad 0 \, , \qquad \qquad \text{otherwise} \, ,
\end{matrix} \right.
\end{equation}
where $\varpi$  is the volume form
on the frame bundle
\begin{equation*} \label{vol}
\varpi \, = \, \bigwedge_{k=1}^n \t^k   \wedge \bigwedge_{(i, j)} \om^i_j  \qquad
\text{(lexicographically ordered)} \, .
\end{equation*}
The invariance property is relative to the \textit{modular character}
 $\d :\Hc_n \ra \Cb$, that extends the \textit{trace} character of
 ${\Fg \Fl} (n, \Rb)$, and is
 defined on generators as follows:
 \begin{equation}
\d(Y_i^j) = \d_i^j , \quad  \d(X_k) = 0, \quad \d(\d_{jk}^i ) =0 ,
 \qquad i, j, k =1, \ldots , n \, .
 \end{equation}
 \medskip
 
\begin{proposition} \label{ibp} 
$1^0$. For any $a, b \in \Ac_n$ and $\, h \in \Hc_n$  one has
\begin{equation} \label{it}
\tau (a\, b) \, = \,  \tau(b \, a) \, , \qquad
\tau (h(a)) \, = \,  \d(h)\, \tau(a) \, .
\end{equation}

$2^0$. For any $h \in \Hc_n$ and $\, a, b  \in \Ac_n$ one has
\begin{equation} \label{sit}
\tau (h(a) \, b) \, = \, \tau \, (a \, \wt{S} (h) (b)) \, , 
\end{equation}
with
\begin{equation} \label{tap}
 \wt{S} (h) \, = \, \sum_{(h)} \, \d(h_{(1)}) \, S(h_{(2)}) \, ,
 \end{equation}
satisfying the anti-involutive property
\begin{equation} \label{ainv}
 \wt{S}^2 \, = \,Id \, .
 \end{equation}
\end{proposition}
\medskip

\begin{proof}  $1^0$. The trace property is a consequence of the 
$\Gc_{n}$-invariance of the volume form $\varpi$. In turn, the 
latter follows from the fact that  
\begin{equation*}
\wt{\vp}^*(\t) = \t \quad \text{and} \quad  \wt{\vp}^* (\om) \, = \, \om \, + \, \g \cdot \t \, ,
\quad \text{(cf. (\ref{omdisp})}\, ; 
 \end{equation*} 
 therefore,
\begin{equation*}
\wt{\vp}^*(\varpi) \, =\,  \bigwedge_{k=1}^n \t^k 
 \wedge \bigwedge_{(i, j)} \left(\om^i_j  + \g^i_{j\ell}(\vp) \t^\ell \right) 
\,  = \, \bigwedge_{k=1}^n \t_k   \wedge \bigwedge_{(i, j)} \om^i_j  \, . 
 \end{equation*}
Passing to the $\Hc_n$-invariance property in (\ref{it}), it suffices to verify it on
generators. Evidently, both sides vanish if $h = \d^i_{jk}$. On the other hand, its
restriction to $\Rb^n \ltimes {\Fg \Fl} (n, \Rb)$ is just the restatement, at the level of
the Lie algebra, of  the right semi-invariance property for the left Haar measure 
on $\Rb^n \ltimes GL (n, \Rb)$.  
\medskip

$2^0$. Using  the `product rule'  (\ref{Yrule}) for vertical vector fields,
in conjunction with
the invariance property (\ref{it}) applied to the product
of two elements $a, b \in \Ac_n$, one obtains
\begin{equation}  \label{Ytap}
\displaystyle \tau \, (Y_i^j (a) \, b ) \, = \,
- \, \tau \, (a \, Y_i^j (b)) \,
 + \, \d_{i}^{j} \, \tau \, (a \, b ) \, , \qquad \fl \, a, b \in \Ac_n \, . \hfill
\end{equation}
On the other hand, 
for the basic horizontal vector fields, (\ref{Xrule}) 
and  (\ref{sit}) give
\begin{equation}  \label{Xtap}
\begin{matrix}
 \displaystyle \tau \, (X_k (a) \, b ) \,  &=& \, - \,
\tau \, (a \, X_k (b) ) \, - \, 
\tau \, ( \d_{jk}^i (a) \, Y_i^j (b) ) \, \hfill \cr 
\displaystyle &=& \, - \, \tau \, (a \, X_k (b) ) \, + \, 
\tau \, (a \, \d_{jk}^i (Y_i^j (b) ) \, ; \hfill  \cr
\end{matrix}
\end{equation}
the second equality uses the $1$-cocycle nature of $\g^i_{jk}$.
The same property implies
\begin{equation}   \label{dtap}
\displaystyle \tau \, (\d^i_{jk} (a) \, b) \, = \,
- \, \tau \, (a \, \d^i_{jk}  (b)) \, , \qquad \fl \, a, b \in \Ac_n
   \, . \hfill
\end{equation}
Thus, the generators of $\Hc_n$ satisfy an \textit{integration by parts} identity
of the form (\ref{sit}). Being multiplicative, this 
rule extends to all elements $h \in \Hc_n$.  Furthermore,
since the pairing $(a, b) \mapsto \tau(a\, b)$ is obviously non-degenerate,
the integration by parts formula uniquely determines an
anti-involutive algebra homomorphism ${\wt S} : \Hc_n \ra \Hc_n$. 
The equations (\ref{Ytap})--(\ref{dtap})
show that  ${\wt S}$ fulfills  (\ref{tap}) on generators, and therefore
for all  $h \in \Hc_n$. 
 \end{proof}
\bigskip

We now define an  \textit{elementary characteristic  cochain}  as a cochain
 $\phi \in C^p (\Ac_n)$  of the form
 \begin{equation*} 
  \phi (a^{0}, \ldots , a^{q})  = \tau
  \left( h^0( a^0) \cdots h^p( a^p)) \right)
   , \, \, h^0 ,..., h^p \in {\Hc}_{n}  ,
   \, \, a^0 ,..., a^p  \in {\Ac}_{n} .
 \end{equation*}
The subspace of $ C^{p}(\Ac_n) $ spanned by
such cochains will be denoted $ C_{\tau}^{p}(\Ac_n) $.
The collection of all characteristic cochains
\begin{equation} \label{cc}
 \Ac_n^{\natural} = \bigoplus_{p \geq 0} C_{\tau}^{p}(\Ac_n)
  \end{equation} 
 forms a module over the cyclic category $\Lb$, more precisely 
 a $\Lb$-submodule of the cyclic module $ \Ac_n^{\natural}$.
Indeed, it is very easy to check that the canonical $\Lb$-operators
\begin{eqnarray} \label{LA} \nonumber
(\d_i \vp) (a^0 , \ldots , a^p) &=& \vp (a^0 , \ldots , a^i \, a^{i+1} , 
\ldots , a^p)\, , \quad i=0,1,\ldots , p-1 \hfill \\  \nonumber
(\d_p  \vp) (a^0 , \ldots , a^p) &=& \vp (a^p \, a^0 , a^1 , \ldots , 
a^{p-1}) \hfill   \\ \nonumber
 (\s_j  \vp) (a^0 , \ldots , a^p) &=& \vp (a^0 , \ldots , a^j , 1 , a^{j+1} 
, \ldots , a^p)\, , \quad j=0,1,\ldots , p \hfill \\ 
(\tau_p \, \vp) (a^0 , \ldots , a^p) &=& \vp (a^p , a^0 , \ldots , a^{p-1}) 
\, . \hfill    
\end{eqnarray}
preserve the property of a cochain of being
characteristic.  
\bigskip

Note that, because of the $\Hc_n$-invariance of the trace (\ref{it}),
$$
C_{\tau}^0 (\Ac_n) \, = \, \Cb \, \tau \, .
$$
More generally, using the integration by parts property
(\ref{sit}), any characteristic cochain can be put in the \textit{standard} form
\begin{equation*} 
  \phi (a^{0}, \ldots , a^{p}) \, = \, 
    \sum_{i=1}^r \, \tau \, (a^{0} \,  h_i^1 (a^1 ) \cdots h_i^p ( a^p ))
  \, , \quad h_{i}^{j} \, \in {\Hc}_{n} \, ;
 \end{equation*} 
moreover, because of the non-degeneracy
of the bilinear form $(a, b)  \mapsto \tau(a\, b) $,
the $p$-differential operator 
$$ 
P (a^{1}, \ldots , a^{p}) \, =  \,
\sum_{i=1}^r \, h_i^1 (a^1 ) \cdots h_i^p ( a^p )
$$
is uniquely determined. In conjunction with  Proposition \ref{multid}, this
means that the linear map
$\, T^{\natural} \, : \, \Hc_n^{ \ot ^{p}}  \longra \, C_{\tau}^{p}(\Ac_n)$,
\begin{equation} \label{Tnat}
T^{\natural} (h^1 \ot \ldots \ot h^{p}) \, : = \, T(1 \ot h^1 \ot \ldots \ot h^{p})
\end{equation} 
is an isomorphism. We can thus transfer via this isomorphism the cyclic 
structure of (\ref{cc}). 
\medskip

\begin{proposition} \label{LH}
The  operators
\begin{eqnarray} \nonumber
\d_0 (h^1 \ot \ldots \ot h^{p-1}) &=& 1 \ot h^1
\ot \ldots \ot h^{p-1} \\ \nonumber 
\d_j (h^1 \ot \ldots \ot h^{p-1}) &=& h^1 \ot \ldots \ot \D h^j \ot
\ldots \ot h^{p-1} , \quad 1 \leq j \leq p-1  \\ \nonumber
\d_p (h^1 \ot \ldots \ot h^{p-1}) &=& h^1 \ot \ldots \ot h^{p-1}
\ot 1   \\ \nonumber
 \s_i (h^1 \ot \ldots \ot h^{p+1}) &=& h^1 \ot \ldots \ot \ve
(h^{i+1}) \ot \ldots \ot h^{p+1} , \quad 0 \leq i \leq p   \\  \nonumber
   \tau_p (h^1 \ot \ldots \ot h^p) &=& (\D^{p-1}  \wt S (h^1)) \cdot h^2
\ot \ldots \ot h^p \ot 1   \, , 
\end{eqnarray}
define a $\Lb$-module structure on
$\, \displaystyle  \Hc_n^{\natural} = 
 \Cb \oplus \bigoplus_{p \geq 1} \Hc_n^{\ot^p} $.
  \end{proposition}
\medskip

\begin{proof}  All we need to check is that the isomorphism
$T^{\natural} $ intertwines the
above $\Lb$-operators with those of (\ref{LA}). This is easy, and  
not surprising, for the 
face operators $\, \{\d_i \, , \, 0 \leq i \leq p-1\}\,$ and the degeneracy operators
$\, \{\s_j \, , \, 0 \leq j \leq p\}$, because it amounts to the natural
equivalence between the cosimplicial structure 
of the \textit{coalgebra} $\Hc_n$ and that of the algebra $\Ac_n$. By contrast,
the full Hopf algebraic structure
of $\Hc_n$ is reflected in the expression of the cyclic operator,
which is obtained as follows:
\begin{eqnarray} \nonumber
\tau_p (T^{\natural} (h^1 \ot \ldots \ot h^{p})) (a^0 , \ldots , a^p) &=&
 T^{\natural} (h^1 \ot \ldots \ot h^{p})(a^p , a^0 , \ldots , a^{p-1}) \\ \nonumber
= \tau (a^{p} \, h^1 (a^0 ) \cdots h^p ( a^{p-1} ))&=&
 \tau (h^1 (a^0 ) \cdots h^p ( a^{p-1} ) \, a^{p})  \\ \nonumber
 &=&  \tau  (a^0 \, \wt S (h^1) \left((h^1)(a^0 ) \cdots h^p ( a^{p-1} ) \, a^{p} \right) ) \, .
 \end{eqnarray}
In the last two lines we have used the trace property and the integration by parts
formula (\ref{tap}). To arrive at the expression in the statement, it remains
to employ the compatibility of $\Hc_n$-action with the
algebra product.
\end{proof}  
\bigskip  

\section{Cyclic cohomology for Hopf algebras}
 
We now consider an arbitrary Hopf algebra $\Hc$ over a field $k$ containing $\Qb$,
with unit $\eta : k \ra \Hc$, counit $\ve : \Hc \ra k$ 
and antipode $S : \Hc \ra \Hc$.  As the discussion of $\Hc_n$ clearly indicates,
the initial datum should also include a `modular' component, playing the role of
the modular function of a Lie group. To be consistent with the intrinsic
duality of the Hopf algebra framework,  this modular datum should also be 
self-dual. One is thus led to postulate the existence of a 
\textit{modular pair}   $(\delta, \sigma)$, consisting of
a character $\delta \in {\Hc}^{*}$,
\begin{equation*}
\delta (a b) = \delta (a) \delta (b) \, , \quad \fl \, a, b \in 
\Hc , \qquad \d(1) = 1 \, ,
\end{equation*}
and a group-like element $\sigma \in \Hc$,
\begin{equation*}
\D(\sigma) = \s \ot \s \, , \qquad \ve (\s) = 1 ,
\end{equation*}
related by the condition
\begin{equation}
    \delta (\s) = 1 .
\end{equation}    
\smallskip

\noindent The convolution of the antipode by the character $\delta$ 
gives rise to the \textit{twisted antipode} $\wt S = S_{\d} : \Hc \ra \Hc $,
\begin{equation}\nonumber
\wt S (h) = \sum_{(h)} \delta (h_{(1)}) \ S (h_{(2)}) \quad , 
\quad h \in \Hc.
\end{equation}
The basic properties of the ordinary antipode are inherited by  the
twisted antipode. In particular, we shall freely use  
the following facts that can be easily checked: 
\begin{eqnarray} \nonumber
\wt S (h^1  h^2) &=& \wt S (h^2)  \wt S (h^1) \, , \quad  \qquad \fl  
h^1 , h^2 \in \Hc \, , \qquad \wt S (1) = 1 \, ; \\ \nonumber
\D  \wt S (h) &=& \sum_{(h)} S (h_{(2)}) \ot \wt S (h_{(1)}) \, , 
\quad \fl  h \in \Hc \, ; \\ \nonumber
\ve \circ \wt S &=& \delta, \quad \quad 
\delta \circ \wt S = \ve. \nonumber
\end{eqnarray}
\smallskip

Since the concrete Hopf algebras often arise as generalized symmetries acting
on algebras, it is reasonable to take into account this dynamical aspect in the
process of defining the  Hopf-cyclic cohomology. 
\medskip

\begin{definition} Let $\Ac$ be an $\Hc$-module algebra.
A linear form $\tau : \Ac \ra k$ will be called a $\s$-{\em trace}
if  
 \begin{equation} \label{str}
\tau (ab) = \tau (b \sigma (a)) \, , \qquad \fl  a,b  \in A  .
\end{equation}
A $\s$-trace $\tau$ will be called $(\Hc, \d)$-{\em invariant} if
\begin{equation} \label{dinv}
\tau (h(a)) \, = \, \d(h) \, \tau (a) \, , \qquad \fl  a  \in A  , \, h
\in  \Hc \, .
\end{equation}
\end{definition}
\medskip

We remark that the invariance condition (\ref{dinv}) is equivalent to the
apparently stronger `integration by parts' property  
\begin{equation} \label{sinv}
\tau (h(a) \, b) \, = \, \tau (a \, \wt S (h) (b)) \, , \qquad \fl  a  \in A  , \, h
\in  \Hc \, .
\end{equation}
Indeed, the former is obtained from the latter by specializing $b=1$.
Conversely, one has (omitting the summation sign in the Sweedler notation)
\begin{eqnarray} \nonumber
\tau (a \, \wt S (h)(b)) &=& \d(h_{(1)}) \, \tau (a \, S(h_{(2)}(b))
\, =  \, \tau \left(\, h_{(1)} (a \,  S(h_{(2)})(b)) \, \right) \\ \nonumber
&=&  \tau \left(h_{(1)} (a )\, h_{(2)} (S(h_{(3)})(b))\right) =
 \tau (h_{(1)} (a )\, \ve(h_{(2)}) \, b) = \tau (h(a) \, b) \, .\nonumber
 \end{eqnarray}
\medskip

Motivated by the discussion in the previous section, we shall enforce
the existence of a natural characteristic homomorphism, as a
guiding principle for the definition of Hopf-cyclic cohomology.
Specifically, the definition should satisfy the following:
\bigskip
 
\textbf{Ansatz.} \textit{Let $(\Ac, \tau)$ be an an $\Hc$-module algebra with
$(\Hc, \d)$-invariant $\sigma$-trace. Then the assignment}
\begin{eqnarray} \nonumber
h^1 \ot \ldots \ot h^n \longmapsto \chi_{\tau} (h^1 \ot \ldots \ot h^n) \in C^n (\Ac)
\, , \quad \fl \, h^1 , \ldots , h^n \in \Hc,\\ \nonumber
\chi_{\tau} (h^1 \ot \ldots \ot h^n) (a^0 , \ldots , a^n) = 
 \tau (a^0  h^1 (a^1) \ldots h^n (a^n)) , \quad 
\fl \, a^0 , \ldots , a^n \in \Ac ,
\end{eqnarray}
\textit{defines a canonical map} $\chi^*_{\tau}: HC^* (\Hc) \ra HC^* (\Ac)$.
 \bigskip
 
If $\Hc$ happens to admit a multi-faithful action on a module algebra, then
the Ansatz would dictate, as  in Proposition \ref{LH},
the following cyclic structure on
\begin{eqnarray} \label{Ans}
 \Hc_{(\delta,\sigma)}^{\natural} \,  &=& \,
 \Cb \oplus \bigoplus_{n \geq 1} \Hc^{\ot^n} \, : \\ \nonumber
\d_0 (h^1 \ot \ldots \ot h^{n-1}) &=& 1 \ot h^1
\ot \ldots \ot h^{n-1} , \\  \nonumber
\d_j (h^1 \ot \ldots \ot h^{n-1}) &=& h^1 \ot \ldots \ot \D h^j \ot
\ldots \ot h^{n-1} , \qquad 1 \leq j \leq n-1  \\  \nonumber
\d_n (h^1 \ot \ldots \ot h^{n-1}) &=& h^1 \ot \ldots \ot h^{n-1}
\ot \s, \\ \nonumber
 \s_i (h^1 \ot \ldots \ot h^{n+1}) &=& h^1 \ot \ldots \ot \ve
(h^{i+1}) \ot \ldots \ot h^{n+1} , \qquad 0 \leq i \leq n \, , \\  \nonumber
   \tau_n (h^1 \ot \ldots \ot h^n) &=& (\D^{p-1}  \wt S (h^1)) \cdot h^2
\ot \ldots \ot h^n \ot \s . 
\end{eqnarray}
The insertion of $\s$ in $\d_n$ and $\tau_n$ accounts for the
passage from an ordinary trace to a twisted trace. 
Note that  
\begin{equation} \nonumber
 {\tau_1}^2 (h) \, = \, {\wt S} ({\wt S}(h) \, \s) \s \, = \, \s^{-1} \,  {\wt S}^2 (h) \, \s ,
 \quad \fl \, h \in \Hc \, ,
\end{equation}
therefore the fact that $\tau_1$ is cyclic is equivalent to the identity
\begin{equation} \label{inv}
 {\wt S}^2  \, = \,  \Ad \s   \, .
\end{equation}

\begin{definition} A modular pair $(\d, \s)$ is said to be {\em in involution}
if the $\d$-twisted antipode  ${\wt S} = S_{\d} $ satisfies the condition
(\ref{inv}).
\end{definition}
\smallskip

The pleasant surprise is that, besides  (\ref{inv}), there are no more
conditions needed for
the existence of a cyclic structure based on the above operators.
The key calculation
to prove this makes the object of the following lemma.
\smallskip

\begin{lemma} \label{cyc}
Given a Hopf algebra $\Hc$ endowed with a modular pair
${(\delta,\sigma)}$ one has, for any $\, h^1, \ldots , h^n \in \Hc$,
\begin{equation} \nonumber
\tau_n^{n+1} (h^1 \ot \ldots \ot h^n) \, = \,
\s^{-1}\, \wt S^2 (h^{1}) \, \s \ot \ldots  \ot \s^{-1}\, \wt S^2 (h^{n}) \, \s  .
\end{equation}
\end{lemma}
\smallskip

\begin{proof}  Omitting the summation sign, we write $\tau_n$ in the form
\begin{equation*}
\tau_n (h^1 \ot h^2 \ot \ldots \ot h^n)  
 = S (h_{(n)}^1)  h^2 \ot S (h_{(n-1)}^1)  h^3 \ot \ldots \ot 
  \wt S (h_{(1)}^1) \s  .
\end{equation*}
Upon iterating once, one obtains $\qquad  \tau_n^2 (h^1 \ot \ldots \ot h^n) \, =$
\begin{eqnarray*}
&=& S (S(h_{(n)}^1)_{(n)}  
h_{(n)}^2 ))  S (h_{(n-1)}^1)  h^3 \ot \cdots 
 \ot \wt S (S (h_{(n)}^1)_{(1)}  h_{(1)}^2) \s \cr \cr
&= & S (h_{(n)}^2)  S (S (h_{(n)(1)}^1))  S (h_{(n-1)}^1) 
 h^3 \ot  \ldots
 \ot \ \wt S (h_{(1)}^2)  \wt S (S(h_{(n)(n)}^1)) \s \cr \cr
&= & S (h_{(n)}^2)  S (h_{(n-1)}^1  S (h_{(n)}^1))  h^3 
\ot \ldots \ot
 \wt S (h_{(1)}^2)  \wt S ( S (h_{(2n-1)}^1 )) \s ;
\end{eqnarray*}
successive use of the basic identities for the counit and the antipode
leads to 
\begin{eqnarray*}
&=& S (h_{(n)}^2)  h^3 \ldots \ot
S (h_{(2)}^2) \ S(S(h_{(2)}^1))
 \wt S (h_{(1)}^1) \s \ot \wt S (h_{(1)}^2)  \wt S ( S (h_{(3)}^1)) \s \cr \cr
&= & \ \D^{(n-1)} \wt S (h^2) \cdot  h^3 \ot  \ldots \ot \s  \ot 
\wt S^2 (h^1) \, \s \,  ;
\end{eqnarray*}
for the last equality we have used, with $k = h_{(1)}^1$,
\begin{eqnarray*} 
S (S(k_{(2)}))  \wt S (k_{(1)}) &=& S (S(k_{(3)})) 
 \d (k_{(1)})  S (k_{(2)}) 
= \d (k_{(1)})  S \left(  k_{(2)}  S (k_{(3)}) 
\right) \cr
&=&  \delta (k_{(1)})  S (\ve (k_{(2)})  1) 
= \delta (k_{(1)})  \ve (k_{(2)}) \cr
&=& \delta \left( k_{(1)}  \ve (k_{(2)}) \right) 
=  \delta (k)  \, .
\end{eqnarray*}

By induction, one obtains 
for any $j = 1, \ldots , n$,
$$
\tau_n^j (h^1 \ot \ldots \ot h^n) = \D^{n-1} \wt S (h^j) \cdot h^{j+1} 
\ot \ldots \ot h^n \ot \s \ot \wt S^2 (h^{1}) \s \ot \ldots \ot \wt S^2 (h^{j-1}) \s.
$$
A last iteration finally gives
$$
\tau_n^{n+1} (h^1 \ot \ldots \ot h^n) = \D^{n-1} \wt S (\sigma) \cdot \wt S^2 
(h^1) \sigma \ot 
\ldots \ot \wt S^2 (h^n) \sigma .
$$
\end{proof}

\medskip

\begin{theorem}
Let $\Hc$ be a Hopf algebra endowed with a modular pair
${(\delta,\sigma)}$.  Then $\,  \Hc_{(\delta,\sigma)}^{\natural} \,  = \,
 \Cb \oplus \bigoplus_{n \geq 1} \Hc^{\ot^n} $
equipped  with the operators (\ref{Ans}) is a $\Lb$-module  
if and only if the modular pair ${(\delta,\sigma)}$ is in involution.
\end{theorem}
\smallskip

\begin{proof} One needs to check that  the operators defined in (\ref{Ans})
satisfy the following relations, that give the standard
presentation of the cyclic category  $\Lb=\D C$ (\cite{Cext}, \cite{L}): 
\begin{equation}\label{ds}
\delta_j  \delta_i = \delta_i  \delta_{j-1}, 
\, \, i < j  , \qquad \s_j  \s_i = 
\s_i  \s_{j+1},  \, \,  i \leq j 
\end{equation}
\begin{equation} \label{sd}
\s_j  \delta_i = \left\{ \begin{matrix}
\delta_i  \s_{j-1} \hfill &i < j \hfill \cr
1_n \hfill &\hbox{if} \ i=j \ \hbox{or} \ i = j+1 \cr
\delta_{i-1}  \s_j \hfill &i > j+1  ;  \hfill \cr
\end{matrix} \right.
\end{equation}
\begin{eqnarray} \label{ci}
\tau_n  \delta_i  = \delta_{i-1}  \tau_{n-1} ,
 \quad && 1 \leq i \leq n ,  \quad \tau_n  \delta_0 = 
\delta_n \\ \label{cj}
\tau_n  \s_i = \s_{i-1} \tau_{n+1} , 
\quad && 1 \leq i \leq n , \quad \tau_n  \s_0 = 
\s_n  \tau_{n+1}^2 \\ \label{ce}
\tau_n^{n+1} &=& 1_n  \, .  
\end{eqnarray}
The cyclicity condition (\ref{ce}) holds by virtue of
Lemma \ref{cyc} and the hypothesis that the modular pair 
${(\delta,\sigma)}$ is in involution. 

The verification of the remaining relations 
is straightforward. As typical illustrations, we check below some
of the compatibility relations. 
For $i=1$ in (\ref{ci}), one has
\begin{eqnarray*}
 &\ &\tau_n \delta_1 (h^1 \ot \ldots \ot h^{n-1}) \, =  \,
 \tau_n (h_{(1)}^1 \ot h_{(2)}^1 \ot h^2 \ot \ldots \ot 
h^{n-1}) \cr
&= & \D^{n-1}  \wt S (h_{(1)}^1) \cdot h_{(2)}^1 \ot h^2 \ot 
\ldots \ot h^{n-1} \ot \s \cr
&= &S (h_{(1)(n)}^1)  h_{(2)}^1 \ot S (h_{(1)(n-1)}^1)  
h^2 \ot \ldots 
\ot S (h_{(1)(2)}^1)  h^{n-1} \ot \wt S (h_{(1)(1)}^1)  \s \cr
&= & \ve (h_{(n)}^1)  1 \ot S (h_{(n-1)}^1)  h^2 \ot \ldots 
\ot S (h_{(1)}^1)  h^{n-1} \ot \wt S (h_{(1)}^1)  \s \cr
&= & \ 1 \ot S (h_{(n-1)}^1)  h^2 \ot \ldots \ot S (h_{(1)}^1)  
h^{n-1} \ot \wt S (h_{(1)}^1) \s \cr
&= & \ \delta_0  \tau_{n-1}  (h^1 \ot \ldots \ot h^{n-1})  , 
\end{eqnarray*}
while for $i=0$ the verification is even easier:
\begin{eqnarray*}
&\ &\tau_n  \delta_0 (h^1 \ot \ldots \ot h^{n-1})\, = \, \ \tau_n (1 \ot 
h^1 \ot \ldots \ot h^{n-1}) = \cr
&= & \ \D^{n-1}  \wt S (1) \cdot h^1 \ot \ldots \ot h^{n-1} \ot \sigma 
=  h^1 \ot \ldots \ot h^{n-1} \ot \sigma \cr
&= & \ \delta_n (h^1 \ot \ldots \ot h^{n-1})  .
\end{eqnarray*}

Passing now to the relations (\ref{cj}), when $i=0$ the left hand side
gives
\begin{eqnarray*}
&&\tau_n  \s_0 (h^1 \ot \ldots \ot h^{n+1}) = \ve (h^1)  \tau_n 
(h^2 \ot \ldots \ot h^{n+1}) = \cr
&&\ =\ \ve (h^1) \sum S (h_{(n)}^2)  h^3 \ot \ldots \ot S (h_{(2)}^2) 
 h^{n+1} \ot \wt S (h_{(1)}^2) \sigma  , 
\end{eqnarray*}
which is reproduced by the right hand side as follows:
\begin{eqnarray*}
& &\s_n  \tau_{n+1}^2 (h^1 \ot \ldots \ot h^{n+1}) = \cr
&= & \ \s_n \left( \sum S (h_{(n+1)}^2)  h^3 \ot \ldots \ot S 
(h_{(2)}^2) \sigma \ot \wt S (h_{(1)}^2)  \wt S^2 (h^1) \sigma \right) \cr
&= & \ \sum \ve (\wt S (h_{(1)}^2) \wt S^2 (h^1) \sigma )  S (h_{(n+1)}^2)  h^3 
\ot \ldots \ot S (h_{(2)}^2) \sigma \cr
&= & \ \ve (\sigma^{-1} \wt S^2 (h^1)\sigma) \sum \delta (h_{(1)}^2)  S (h_{(n+1)}^2)  h^3 \ot 
\ldots \ot S (h_{(2)}^2) \sigma \cr
&= & \ \ve (h^1)  S (h_{(n)}^2)  h^3 \ot \ldots \ot S (h_{(2)}^2) 
 h^{n+1} \ot \wt S (h_{(1)}^2) \sigma . 
\end{eqnarray*}
For $i=1$ the left hand side is
\begin{eqnarray*}
& &\ \tau_n  \s_1 (h^1 \ot \ldots \ot h^{n+1}) =  \ve (h^2)  
\tau_n (h^1 \ot h^3 \ot \ldots \ot h^{n+1}) \cr
&= & \ \ve (h^2) \cdot \D^{n-1}  \wt S (h^1) \cdot h^3 \ot \ldots 
\ot h^{n+1} \ot \sigma  , 
\end{eqnarray*}
and the right hand side amounts to the same:
\begin{eqnarray*}
&& \ \s_0  \tau_{n+1} (h^1 \ot \ldots \ot h^{n+1}) = \cr
&& \ \sum \s_0 (S (h_{(n+1)}^1)  h^2 \ot \ldots \ot S (h_{(2)}^1)  
h^{n+1} \ot \wt S (h_{(1)}^1)) \sigma \cr
&= & \ \sum \ve (h^2) \cdot \ve (h_{(n+1)}^1) \cdot S (h_{(n)}^1)  h^3 
\ot \ldots \ot S (h_{(2)}^1)  h^{n+1} \ot \wt S (h_{(1)}^1) \sigma \cr
&= & \ \sum \ve (h^2) \cdot S (h_{(n-1)}^1)  h^3 \ot \ldots \ot S 
(h_{(2)}^1)  h^{n+1} \ot \wt S (h_{(1)}^1) \sigma  .
\end{eqnarray*}
Similar calculations hold for $i = 2, \ldots n$.
\end{proof}
\medskip

For completeness, we record below the normalized
bi-complex  
$$\, (C C^{*, *} (\Hc ; \delta,\sigma), \, b , \, B )$$
that computes the Hopf-cyclic cohomology of $\Hc$
with respect to a modular pair in involution ${(\delta,\sigma)}$:
\begin{equation*} \label{CbB}
CC^{p, q} (\Hc; \delta,\sigma)  \, = \, \left\{
\begin{matrix} {\bar C}^{q-p} (\Hc; \delta,\sigma) \, , \quad q \geq p \, , \cr
  0 \, ,  \quad \qquad  \qquad q <  p \, , \end{matrix}   \right.
\end{equation*}
where
\begin{equation*} 
 {\bar C}^{n} (\Hc; \delta,\sigma) \, = \, \left\{
\begin{matrix} \bigcap_{i=0}^{n-1} \Ker \,\s_i  \, , \quad n \geq 1 \, , \cr \cr
  \Cb \, ,  \quad \qquad  \qquad n = 0 \, ; \end{matrix}   \right.
\end{equation*}
 the operator
\begin{equation*}  \nonumber
b: {\bar C}^{n-1} (\Hc; \delta,\sigma) \ra {\bar C}^n (\Hc; \delta,\sigma), 
\qquad b =
\sum_{i=0}^{n} (-1)^i \d_i \,
\end{equation*}
has the form $\, \qquad b ( \Cb ) \, = \, 0 \, \,$ for $\, n = 0 \,$, 
\begin{eqnarray} \nonumber
\ b\, (h^1 \ot \ldots \ot h^{n-1}) &=& 1 \ot h^1 \ot \ldots \ot
h^{n-1}  \cr & + & \ \sum_{j=1}^{n-1} (-1)^j \sum_{(h_{j})} h^1
\ot \ldots \ot h_{(1)}^j \ot h_{(2)}^j \ot \dots \ot h^{n-1} \cr
&+ & \ (-1)^n h^1 \ot \ldots \ot h^{n-1} \ot 1 \, ,
\end{eqnarray}
while the $B$-operator $B: {\bar C}^{n+1} (\Hc; \delta,\sigma) \ra 
{\bar C}^n (\Hc; \delta,\sigma)$ 
is defined by the formula
\begin{equation} \nonumber
  B = A \circ B_{0} \, ,  \quad n \geq 0 \, ,
\end{equation}
where  
\begin{equation*} 
B_{0} (h^1 \ot \ldots \ot h^{n+1})  \, = \,\left\{
 \begin{matrix} (\D^{n-1}
  \wt S (h^1)) \cdot h^2
\ot \ldots \ot h^{n+1} \, , \qquad   n \geq 1 \, , \cr
\d (h^1) \, , \qquad \qquad \qquad \qquad \qquad \text{for} \quad n = 0  
\end{matrix} \right.
\end{equation*}
 and
\begin{equation} \nonumber
A = 1 + \lb_n  + \cdots +  \lb_n^n \,  , \qquad
\text{with} \qquad \lb_n = (-1)^n \tau_n \, .
\end{equation}

The groups $\, \{H C^n  (\Hc; \delta,\sigma)\}_{n \in \Nb} \,$ are computed from the
first quadrant total complex $\, ( TC^{*} (\Hc; \delta,\sigma),\, b+B ) \,$,
\begin{equation*}
    TC^{n}(\Hc; \delta,\sigma) \,  = \, \sum_{k=0}^{n} \, C C^{k, n-k} (\Hc; \delta,\sigma) \, ,
\end{equation*}
and the periodic groups $\, \{H P^i  (\Hc; \delta,\sigma)\}_{i \in \Zb/2} \,$ 
are computed from the full total complex 
$\,(T P^{*}(\Hc; \delta,\sigma),\,  b+B ) \,$,
\begin{equation*}
      T P^{i}(\Hc; \delta,\sigma) \,  = \, \sum_{k \in \Zb} \, C C^{k, i-k} (\Hc; \delta,\sigma) \, .
\end{equation*}
\bigskip

A  class of useful examples is that of 
universal enveloping algebras of Lie algebras. 
If $\Hc = {\FA}({\Fg})$,
with ${\Fg}$ a Lie algebra over $k$, together with a character 
$\delta : \Fg \ra k$, then
the periodic Hopf-cyclic cohomology of the
$\left({\FA}({\Fg}); \d, 1\right)$ reduces to the Lie algebra
homology with coefficients in the 
$1$-dimensional $\FG$-module ${k}_{\d}$ associated
to the character (cf. \cite[$\S 7$, Prop. 7]{CMhti}, see also
Thm. \ref{Lie} below):
\begin{equation} \nonumber
 HP^* \left({\FA}({\Fg}) ; \d, 1 \right) \simeq 
 \sum^{\oplus}_{i \equiv * \, (2)} \quad H_i \, ({\Fg} , {k}_{\d}) .
\end{equation}
Dually, if $\Hc = {\Hc} (G)$ is
the Hopf algebra of regular (i.e. polynomial) functions
on a unipotent affine algebraic group $G$
with Lie algebra ${\Fg}$,
then its periodic cyclic cohomology with respect to the trivial
modular pair $(\ve, 1)$ is isomorphic to the Lie algebra cohomology:
\begin{equation} \nonumber
  HP^* \left({\Hc} (G) ; \ve, 1\right) \simeq  
 \sum^{\oplus}_{i \equiv * \, (2)} \quad H^{i} \, ({\Fg}, k) .
\end{equation}
\smallskip

In these instances $\s = 1$. 
On the other hand, the class of
quasitriangular Hopf algebras called   
\textit{ribbon algebras}, and in particular the $q$-deformations 
${\FA}_q({\Fg})$ of the
classical Lie algebras, provide examples of dual pairs
with $\s \neq 1$ and $\d = \ve$ (see \cite[\S 4]{CMcchs}).
\medskip

Our prototype of a Hopf algebra with
modular pair in involution was of course $(\Hc_n ; \d, 1)$.  
In \cite[$\S 7$, Theorem 11]{CMhti} we established
a canonical isomorphism
\begin{equation}  \label{gf}
\kappa_{n}^{*} \, : \sum^{\oplus}_{i \equiv * \, (2)} \quad  H^{i} ({\Fa}_n ; \mathbb{C})
\build\longrightarrow_{}^{\simeq} HP^{*} \, (\Hc_n ; \d, 1) \, ,
\end{equation}
between the Gelfand-Fuks cohomology of the Lie algebra  ${\Fa}_n $ of formal
vector fields on $\mathbb{R}^{n}$ and the periodic Hopf-cyclic cohomology
of $(\Hc_n ; \d, 1)$. Moreover,
the isomorphism $\kappa_{n}^{*} $
is implemented at the cochain level by a
homomorphism constructed out of a (fixed but arbitrary)
torsion-free connection. 
\smallskip

A similar statement takes place at the \textit{relative} level, and it actually plays a
key role in understanding the Chern character of the hypoelliptic signature
operator, namely the existence
of a canonical isomorphim
\begin{equation}  \label{rgf}  
\kappa_{n, \SO(n)}^{*} \, : \sum^{\oplus}_{i \equiv * \, (2)} \quad  
H^{i} ({\Fa}_n , \SO(n) ;\mathbb{C})
\build\longrightarrow_{}^{\simeq} HP^{*} \, (\Hc_n , \SO(n) ; \d, 1) \, 
\end{equation}
While the meaning of the relative Hopf cyclic group $HP^{*} \, (\Hc_n , \SO(n) ; \d, 1)$
happened  to be quite clear in that particular context, it was not so in the case of
non-compact isotropy, for instance  for the Lorentz group $\SO(n-1, 1)$. 
The treatment of the general case makes the object of the next section.  
\bigskip

\section{Relative Hopf cyclic cohomology}
 
Let $\Hc$ be an arbitrary Hopf algebra over a field $F$ containing $\Qb$,
with unit $\eta : F \ra \Hc$, counit $\ve : \Hc \ra F$ 
and antipode $S : \Hc \ra \Hc$, and let  $\Kc$ be a Hopf subalgebra of  $\Hc$.  
We consider the tensor product over $\Kc$
\begin{equation}
	\Cc \, = \, \Cc (\Hc, \Kc) \, := \, \Hc \ot_{\Kc} F   \, ,
\end{equation}
where $\Kc$ acts on $\Hc$ by right multiplication and on $F$ by the
counit. It is a left $\Hc$-module, which can be identified with the
quotient module $\, \Hc/\Hc \Kc^+ $, $\, \Kc^+ = \Ker \ve | \Kc \,$, via the
isomorphism  induced by the map 
\begin{equation}
h \in \Hc \, \longmapsto \, \dot{h} = h \ot_{\Kc} 1 \in \Hc \ot_{\Kc} F \, . 
\end{equation}
Moreover, $\, \Cc = \Cc (\Hc, \Kc) \,$ is an $\Hc$-module coalgebra. Indeed,
its coalgebra structure is given by the
coproduct  
 \begin{equation} \label{copr}
	 \D_{\Cc} \, (h \ot_{\Kc} 1) \, = \, (h_{(1)} \ot_{\Kc} 1) \ot (h_{(2)} \ot_{\Kc} 1) \, ,
\end{equation}
inherited from that on $\Hc$,
$\, \,	 \D h  \, = \, h_{(1)} \ot  h_{(2)} \, $ and is compatible 
with the action of
$\Hc$ on $\Cc$ by left multiplication:
\begin{equation*}
 \D_{\Cc} \, (h' \cdot h \ot_{\Kc} 1) \, = \,  \D h'  \cdot  \D_{\Cc} \, (h \ot_{\Kc} 1) \, ;
 \end{equation*}
similarly, there is an inherited counit
\begin{equation} \label{coun}
	 \ve_{\Cc} \, (h \ot_{\Kc} 1) \, = \,  \ve (h) \, , \quad \fl \, h \in \Hc \, ,
\end{equation}
that satisfies
\begin{equation*}
  \ve_{\Cc} \, (h' \cdot h \ot_{\Kc} 1) \, = \,  \ve (h) \, \ve_{\Cc} \, (h \ot_{\Kc} 1) \,  .
 \end{equation*} 
\medskip

Let us now fix a right $\Hc$-module $M$ that is also a left $\Hc$-comodule, 
\begin{equation*}
  {}_M \D (m) \, = \,  m_{(-1)}   \ot m_{(0)}  \, \in \Hc \ot M \,  ,
 \end{equation*} 
and assume that
it is a \textit{stable anti-Yetter-Drinfeld module}, cf. \cite{HKRS1}, that is
\begin{eqnarray} \label{stable}
m_{(0)} \, m_{(-1)} \, &=& \,  m \, ;  \\ \label{ayd}
 {}_M \D (m \, h) \, &=& \,  S(h_{(3)} ) \, m_{(-1)} \, h_{(1)}  \ot  m_{(0)}\,  h_{(2)} \, .
 \end{eqnarray}
\medskip

The following
result provides, in the extended framework of  \cite{HKRS2},
a direct generalization to the relative case with coefficients
of the definition (\ref{Ans}) suggested by the Ansatz of $\S 4$.
\medskip

\begin{theorem} \label{rcc}
Let 
$C^{\ast} (\Hc, \Kc ; M) = \{ C^n  (\Hc, \Kc ; M) 
: = \, M \ot_{\Kc} \Cc^{\ot n} \}_{n \geq 0}$,
where $\Cc = \Hc \ot_{\Kc} F $ and the tensor product over
$\Kc$ is with respect to the diagonal action on $ \Cc^{\ot n}$. 
 The following operators are well-defined and 
endow $C^{\ast} (\Hc, \Kc ; M) $ with a cyclic structure:
 \begin{eqnarray} \label{cmod}
 \d_0 (m \ot_{\Kc} c^1 \ot \ldots \ot c^{n-1}) &=& 
 m \ot_{\Kc} \dot{1} \ot c^1 \ot \ldots \ot
  \ldots \ot c^{n-1} ,   \\  \nonumber
\d_i  (m \ot_{\Kc} c^1 \ot \ldots \ot c^{n-1})  &=& 
m \ot_{\Kc} c^1 \ot \ldots \ot
 c^i_{(1)} \ot c^{i}_{(2)} \ot
\ldots \ot c^{n-1} ,   \\  \nonumber
&\,& \qquad  \qquad  \qquad \qquad  \fl \quad 1 \leq i \leq n-1 \, ; \\  \nonumber
\d_n  (m \ot_{\Kc} c^1 \ot \ldots \ot c^{n-1})  &=& 
m_{(0)} \ot_{\Kc}  c^1 \ot \ldots \ot c^{n-1}
\ot \dot{m}_{(-1)}\, ; \\ \nonumber
 \s_i  (m \ot_{\Kc} c^1 \ot \ldots \ot c^{n+1})  &=& 
 m \ot_{\Kc}  c^1 \ot \ldots \ot \ve
(c^{i+1}) \ot \ldots \ot c^{n+1} ,  \\ \nonumber
&\,& \qquad  \qquad  \qquad\qquad  \fl \quad 0 \leq i \leq n \, ; \\    \nonumber
   \tau_n  (m \ot_{\Kc}  \dot{h}^1 \ot c^2 \ot \ldots \ot c^n) &=& 
   m_{(0)} h^1_{(1)} \ot_{\Kc}  S( h^1_{(2)}) \cdot (c^2
   \ot \ldots \ot c^n \ot \dot{m}_{(-1)}).
\end{eqnarray}
\end{theorem}
\smallskip

\begin{proof} One can directly check that the above operators are defined in a
consistent fashion and satisfy the defining relations of the cyclic category.
This is however redundant, because they can simply be obtained by transport of
structure. Indeed,
since $\, \Cc \,$ is an $\Hc$-module coalgebra, one can apply
 \cite[Theorem 2.1]{HKRS2} that provides a cyclic structure on
 the collection of spaces
\begin{equation} \label{cs}
 \tilde{C}^{\ast} (\Hc, \Kc ; M) = \{ \tilde{C}^n  (\Hc, \Kc ; M) 
: = \, M \ot_{\Hc} \Cc^{\ot n+1} \}_{n \geq 0}  \, \, ,
\end{equation}
where $\Hc$ acts diagonally on $ \Cc^{\ot n+1}$, as follows:
 \begin{eqnarray} \label{ccmc} \nonumber
\tilde{\d}_i (m \ot_{\Hc} c^0 \ot \ldots \ot c^{n-1}) &=& m \ot_{\Hc} c^0 \ot \ldots \ot
 c^i_{(1)} \ot c^{i}_{(2)} \ot
\ldots \ot c^{n-1} ,   \\  \nonumber
&\,& \qquad  \qquad  \qquad \qquad  \fl \quad 0 \leq i \leq n-1 \, ; \\  \nonumber
\tilde{\d}_n (m \ot_{\Hc} c^0 \ot \ldots \ot c^{n-1}) &=& 
m_{(0)} \ot_{\Hc}  c^0_{(2)} \ot c^1 \ot \ldots \ot c^{n-1}
\ot m_{(-1)} c^0_{(1)} ; \\ \nonumber
 \tilde{\s}_i (m \ot_{\Hc} c^0 \ot \ldots \ot c^{n+1}) &=& m \ot_{\Hc}  c^0 \ot \ldots \ot \ve_{\Cc}
(c^{i+1}) \ot \ldots \ot c^{n+1} ,  \\  \nonumber
&\,& \qquad  \qquad  \qquad\qquad  \fl \quad 0 \leq i \leq n \, ; \\  \nonumber
   \tilde{\tau}_n (m \ot_{\Hc} c^0 \ot \ldots \ot c^{n}) &=& 
   m_{(0)} \ot_{\Hc}  c^1 \ot \ldots \ot c^{n} \ot m_{(-1)} c^0  . 
\end{eqnarray}
To recast the cyclic module (\ref{cs}) in the form that appears in the statement
we shall make use of the  `transfer' isomorphism 
\begin{equation*} 
	 \Psi^{\ast} : C^{\ast} (\Hc, \Kc ; M) \longrightarrow  
	 \tilde{C}^{\ast} (\Hc, \Kc ; M) \, 
\end{equation*}
whose $n$-th component is defined by the simple formula
\begin{equation} \label{trans}
	 \Psi^n   (m \ot_{\Kc} c^1 \ot \ldots \ot c^n)   \, = \,
	 m \ot_{\Hc}   \dot{1}  \ot c^1 \ot \ldots \ot c^n \, ,
\end{equation}
which is obviously consistent.
\smallskip

We claim that the expression 
 \begin{equation} \label{transinv}
	 \Phi^n  (m \ot_{\Hc} \dot{h}^0 \ot c^1 \ot \ldots \ot c^n)   \, = \,
	 m\, h^0_{(1)} \ot_{\Kc} S(h^0_{(2)})  \cdot (c^1 \ot \ldots \ot c^n) ,
\end{equation}
 is also well-defined
and gives the inverse operator to $ \Psi^n$. 
Note first that the right hand side only
depends on the class $c^0 = \dot{h}^0 \in \Cc$; indeed, for any 
$g \in \Hc$ and $k \in \Kc^+$,
\begin{eqnarray*} 
m\, g^0_{(1)} k_{(1)}  \ot_{\Kc} S( k_{(2)})
S(g^0_{(2)})  \cdot (c^1 \ot \ldots \ot c^n) \, &=& \\
m\, g^0_{(1)} k_{(1)} S( k_{(2)})  \ot_{\Kc} 
S(g^0_{(2)})  \cdot (c^1 \ot \ldots \ot c^n) \, &=& \\
m\, g^0_{(1)} \ve (k)  \ot_{\Kc} 
S(g^0_{(2)})  \cdot (c^1 \ot \ldots \ot c^n) \, &=&  \, 0 \, .
\end{eqnarray*}
Let us now check the consistency of the definition (\ref{transinv}). On replacing
the elementary tensor $ \dot{h}^0 \ot c^1 \ot \ldots \ot c^n \, \in \,  \Cc^{\ot n+1} $ 
in the left hand side by 
\begin{equation*} 
 h \cdot (\dot{h}^0 \ot c^1 \ot \ldots \ot c^n) \, = \,
\dot{\widehat{h_{(1)} h^0}} \ot  h_{(2)} \cdot ( c^1 \ot \ldots \ot c^n) \, ,
 \qquad h \in \Hc \, ,
\end{equation*} one obtains in the right hand side
\begin{eqnarray*} 
m\,  h_{(1)}  h^0_{(1)} \ot_{\Kc}  S(h^0_{(2)}) S(h_{(2)}) \, h_{(3)}
\cdot (c^1 \ot \ldots \ot c^n )\, &=& \, \\   
m\,  h_{(1)}  h^0_{(1)} \ot_{\Kc}  S(h^0_{(2)}) \ve (h_{(2)})
\cdot (c^1 \ot \ldots \ot c^n ) \, &=& \, \\   
m\,  h_{(1)} \ve (h_{(2)}) h^0_{(1)} \ot_{\Kc}  S(h^0_{(2)}) 
\cdot (c^1 \ot \ldots \ot c^n ) \, &=& \, \\   
m\,  h h^0_{(1)} \ot_{\Kc}  S(h^0_{(2)}) 
\cdot (c^1 \ot \ldots \ot c^n ) &,&
\end{eqnarray*}
which correctly corresponds to $\,  \Phi^n 
( m h \ot_{\Hc} \dot{h}^0 \ot c^1 \ot \ldots \ot c^n )$.
\smallskip

Next, $ \Phi^n$ is obviously a left
inverse to $ \Psi^n$. Conversely, one has 
\begin{eqnarray*} 
( \Psi^n \circ \Phi^n) (m \ot_{\Hc} \dot{h}^0 \ot \ldots \ot \dot{h}^n) &=&
 m\, h^0_{(1)} \ot_{\Hc}   \dot{1}  \ot S(h^0_{(2)})  \cdot (\dot{h}^1 \ot \ldots \ot \dot{h}^n) \\
 &=& m \ot_{\Hc} 
 \dot{\widehat{h^0_{(1)}}} \ot h^0_{(2)} S(h^0_{(3)})  \cdot (\dot{h}^1 \ot \ldots \ot \dot{h}^n) \\
 &=& m \ot_{\Hc} 
\dot{ \widehat{h^0_{(1)}}} \ve (h^0_{(2)}) \ot  \dot{h}^1 \ot \ldots \ot \dot{h}^n \\
&=& m \ot_{\Hc} \dot{h}^0  \ot  \dot{h}^1 \ot \ldots \ot \dot{h}^n \, . 
\end{eqnarray*}

To achieve the proof, it remains to notice  that defining the face, degeneracy
and cyclic operators by transport of structure,
 \begin{eqnarray*} 
  \d_i  &=&\Phi^n \circ \tilde{\d}_i  \circ \Psi^n  \, , \qquad \qquad
   0 \leq i \leq n-1 \, ; \\  \nonumber
 \s_i  &=& \Phi^n \circ \tilde{\s}_i  \circ \Psi^n  \, ,  \qquad 
  \qquad 0 \leq i \leq n \, ; \\  \nonumber 
 \tau_n  &=&   \Phi^n \circ \tilde{\tau}_n \circ \Psi^n \, ,
 \end{eqnarray*}
 their expressions are precisely as stated.
\end{proof}
\bigskip

\begin{definition}
The relative Hopf cyclic cohomology $HC^{\ast} (\Hc, \Kc ; M)$
of the pair $\Kc \subset \Hc$ \textit{with coefficients} in the stable anti-Yetter-Drinfeld 
module-comodule $M$ is the cyclic cohomology of the
cyclic module (\ref{cmod}). 
\end{definition}

Dually, one can define the 
 \textit{relative Hopf cyclic homology} $HC_{\ast} (\Hc, \Kc ; M') $ 
 of the pair $(\Hc, \Kc)$ \textit{with coefficients} in the dual
 module-comodule $\, M' = \Hom (M, F) \, $ as
 the cyclic homology of the cyclic module 
\begin{equation} \label{rcs}
 C_{\ast} (\Hc, \Kc ; M' ) = \{ C_n  (\Hc, \Kc ; M' ) 
: = \, \Hom_{\Kc} (\Cc^{\ot n},  M')  \}_{n \geq 0}  \, ,
\end{equation} 
obtained by dualizing the cyclic structure of Theorem \ref{rcc}:
 \begin{eqnarray*} 
 (d_0 \vp )( c^1 \ot \ldots \ot c^{n-1}) &=& 
 \vp ( \dot{1} \ot c^1 \ot \ldots \ot
  \ldots \ot c^{n-1}),   \\ 
(d_i \vp) (c^1 \ot \ldots \ot c^{n-1})  &=& 
\vp ( c^1 \ot \ldots \ot
 c^i_{(1)} \ot c^{i}_{(2)} \ot
\ldots \ot c^{n-1}) ,   \\ 
&\,& \qquad  \qquad  \qquad \qquad  \fl \quad 1 \leq i \leq n-1 \, ; \\ 
(d_n \vp) (c^1 \ot \ldots \ot c^{n-1}) (m) &=& 
\vp ( c^1 \ot \ldots \ot c^{n-1}
\ot \dot{m}_{(-1)}) (m_{(0)}) \, ; \\ \nonumber
\\
 (s_i \vp)  ( c^1 \ot \ldots \ot c^{n+1})  &=& 
\vp ( c^1 \ot \ldots \ot \ve
(c^{i+1}) \ot \ldots \ot c^{n+1}) ,  \\  \nonumber
&\,& \qquad  \qquad  \qquad\qquad  \fl \quad 0 \leq i \leq n \, ; \\    
  (t_n \vp)  ( \dot{h}^1 \ot c^2 \ot \ldots \ot c^n) (m) &=&
 \vp \left(S( h^1_{(2)}) \cdot (c^2
   \ot \ldots \ot c^n \ot \dot{m}_{(-1)}) \right)(m_{(0)} h^1_{(1)}). 
\end{eqnarray*}
 \bigskip
 
 \begin{remark} \label{adj}
 The restriction to $\Kc$ of the
 left action of $\Hc$ on $\Cc = \Hc \ot_{\Kc} F $ can also be regarded
 as `adjoint action', induced by conjugation.
 \end{remark}

Indeed, for $k \in \Kc$,
\begin{equation*}
k_{(1)}\,  h\,  S(k_{(2)} ) \ot_{\Kc} 1 \, = \, k_{(1)}\, h  \ot_{\Kc} \ve (k_{(2)} )\, 1 
\, = \, k \, h \ot_{\Kc} 1\, .
\end{equation*}

\bigskip

We now specialize the above notions to the Lie algebra case, in order
 to verify that they coincide with the usual definitions of
 relative Lie algebra homology and
cohomology.
\medskip

Let $\Fg$ be a Lie algebra over the field $F$, let $\Fh \sbs \Fg$ be a 
Lie subalgebra,
and let $M $ be a $\Fg$-module.  We equip $M $ with the
\textit{trivial} $\Fg$-comodule structure
\begin{equation} \label{cotriv}
  {}_M \D (m) \, = \,  1   \ot m  \, \in \Hc \ot M \,  ,
 \end{equation} 
and note the condition
(\ref{stable}) is then trivially satisfied, while (\ref{ayd}) follows from
(\ref{cotriv}) and
the cocommutativity of the universal enveloping algebra $\FA (\Fg)$. 
The relative Lie algebra homology and cohomology
of the pair $ \Fh \sbs \Fg$ with coefficients in $M$,
resp. $M'$, 
are computed from the Chevalley-Eilenberg (cf. \cite{CE}) complexes
\begin{eqnarray*} 
\{C_* (\Fg, \Fh ; M) , {}_M\d \} , \quad 
C_n (\Fg, \Fh ; M) &:=& M \ot_{\Fh} \bigwedge ^n (\Fg/\Fh) \, , \\
\text{resp.} \quad \{C^* (\Fg, \Fh ; M') , d_{M'} \} ,
  \quad  C^n (\Fg, \Fh ; M')  &:=& \Hom_{\Fh}
 \left(\bigwedge ^n (\Fg/\Fh) , M' \right) \, ;
\end{eqnarray*}
the action of $\Fh$ on $\Fg/\Fh$ is induced by the adjoint 
representation and the differentials are given by the formulae
\begin{eqnarray} \nonumber
{}_M\d (m\ot_{\Fh}&&\dot{X_1} \wdg \ldots \wdg \dot{X}_{n+1} ) \, = \,
\sum_{i=1}^{n+1}  (-1)^{i + 1}  m X_i \ot_{\Fh}
\dot{X}_1 \wdg \ldots \wdg \Check{\dot{X}}_i \ldots  \wdg \dot{X}_{n+1}  \\ \label{dh} 
+ &&  \sum_{i < j} (-1)^{i + j} m \ot_{\Fh} \dot{\widehat{[X_i , X_j ]}} \wdg \dot{X}_1 
\wdg \ldots \wdg \Check{\dot{X}}_i \ldots \Check{\dot{X}}_j \ldots  \wdg \dot{X}_{n+1} 
\\ \nonumber
(d_{M'}\vp)&& (\dot{X_1} \wdg \ldots \wdg \dot{X}_{n+1} ) \, = \,
\sum_{i=1}^{n+1}  (-1)^{i + 1}  X_i \cdot \vp 
( \dot{X}_1 \wdg \ldots \wdg \Check{\dot{X}}_i \ldots  \wdg \dot{X}_{n+1} ) \\ \label{dc} 
+ && \sum_{i < j} (-1)^{i + j+1} \vp ( \dot{\widehat{[X_i , X_j ]}} \wdg \dot{X}_1 
\wdg \ldots \wdg \Check{\dot{X}}_i \ldots \Check{\dot{X}}_j \ldots  \wdg \dot{X}_{n+1} )
\end{eqnarray}
where $\dot{X} \in \Fg/\Fh$ stands for the class modulo $\Fh$
of $X \in \Fg$ and the superscript 
$\Check{}$ signifies the omission of the indicated variable.
\bigskip

\begin{theorem} \label{Lie}
There are canonical isomorphisms between the periodic relative Hopf cyclic
cohomology (resp. homology) of the pair $\FA(\Fh) \subset \FA(\Fg) $, with
coefficients in any $\Fg$-module $M$ (resp. its dual $M'$),
and the relative Lie algebra homology (resp. cohomology) with coefficients
of the pair
$\Fh \subset \Fg $ :
\begin{eqnarray*}
HP^* (\FA(\Fg), \FA(\Fh) ; M) &\cong& \bigoplus_{n  \equiv * \mod 2} 
H_n (\Fg, \Fh ; M)   \\
\text{resp.} \quad HP_* (\FA(\Fg), \FA(\Fh) ; M') &\cong& \bigoplus_{n  \equiv * \mod 2} 
H^n (\Fg, \Fh ; M') .
\end{eqnarray*} 
\end{theorem}

\begin{proof} The pattern of the proof is exactly the same as in the ``absolute'' case, cf.
\cite[\S 7, Prop.7]{CMhti}. For the convenience of the reader,
we shall repeat the main steps in order to make sure that
they apply in the relative case as well.
\smallskip

First of all, the Hopf cyclic cohomology $HC^* (\FA(\Fg), \FA(\Fh) ; M)$ can be computed
from the normalized
bi-complex $\left( \bar{C}^* (\FA(\Fg), \FA(\Fh) ; M) , {}_M b, {}_M B \right)$ where 
\begin{equation*} 
{\bar C}^n(\FA(\Fg), \FA(\Fh) ; M) \, = \,
 \bigcap_{i=0}^{n-1} \;{\rm Ker}\,\s_i \,, \quad \fl n \geq 1, \qquad 
                        {\bar C}^0 (\FA(\Fg), \FA(\Fh) ; M)  = F ;
    \end{equation*}
the operator $\, {}_M b = \sum_{i=0}^{n} (-1)^i \d_i : {\bar C}^{n-1}  \ra {\bar C}^n $ 
has the expression
\begin{eqnarray} \label{bM} \nonumber
{}_M b (m &\ot_{\Kc}& c^1 \ot \ldots \ot c^{n-1}) \, =\, m \ot_{\Kc} \dot{1} \ot c^1 \ot \ldots \ot
  \ldots \ot c^{n-1} \\  \nonumber
& + & \ \sum_{i=1}^{n-1} (-1)^i m \ot_{\Kc} c^1 \ot \ldots \ot
 c^i_{(1)} \ot c^{i}_{(2)} \ot
\ldots \ot c^{n-1} ,  \\  
&+ & \ (-1)^n m \ot_{\Kc}  c^1 \ot \ldots \ot c^{n-1} \ot \dot{1}
\end{eqnarray}
involving only the coalgebra structure, while
$\, {}_MB: {\bar C}^{n+1}  \ra {\bar C}^n $ is given by the formula
  \begin{equation} \nonumber
 {}_MB = \sum_{i=0}^n (-1)^{ni} \, {\tau_n}^i \circ \s_{-1} \, , \qquad
 \text{where} \qquad \s_{-1}  = \s_n \circ \tau_{n+1} 
 \end{equation}
is the extra degeneracy operator
 \begin{equation} \nonumber
 \s_{-1}(m \ot_{\Kc}  \dot{h}^1 \ot c^2 \ot \ldots \ot c^{n+1}) =
   m h^1_{(1)} \ot_{\Kc}  S( h^1_{(2)}) \cdot (c^2
   \ot \ldots \ot c^{n+1}) .
\end{equation}
\smallskip

We recall that the antisymmetrization map  $\, \a^n : \bigwedge ^n (\Fg/\Fh)  \ra  
\Cc(\FA(\Fg), \FA(\Fh))^{\ot n}$,
\begin{equation*}
\a^n (\dot{X_1} \wdg \ldots \wdg \dot{X}_n ) \, := \, \frac{1}{n!}
\sum_{\s \in S_n} \sign (\s)  \dot{X}_{\s(1)} \ot \ldots \ot \dot{X}_{\s(n)} ,
\end{equation*} 
has a left inverse 
\begin{equation*}
 \mu^n:  \Cc(\FA(\Fg), \FA(\Fh))^{\ot n}  \ra  \bigwedge ^n (\Fg/\Fh) \, , \qquad
 \mu^n \circ \a^n = \Id \, ,
 \end{equation*} 
which is the unique DGA map whose first-degree component  
$\mu^1$ is the canonical projection of $\Cc = \Cc(\FA(\Fg), \FA(\Fh)) \simeq
S(\Fg)/S(\Fg) S(\Fh)^+$ onto $\Fg/\Fh$.

In view of Remark \ref{adj},   
both maps commute with the adjoint representation, and 
therefore induce well-defined maps
\begin{eqnarray*}  
{}_M \a^n : C_n (\Fg, \Fh ; M) = M \ot_{\Fh} \bigwedge ^n (\Fg/\Fh)  &\ra&  
{\bar C}^n  (\FA(\Fg), \FA(\Fh) ; M) \, , \\
\text{resp.} \quad
{}_M \mu^n  : M \ot_{\FA(\Fh)} \Cc^{\ot n} &\ra&  M \ot_{\Fh} \bigwedge ^n (\Fg/\Fh) .
\end{eqnarray*} 
By transposition one also obtains dual maps 
\begin{eqnarray*}  
 \mu^n_{M'} :  \Hom_{\Fh} \left(\bigwedge ^n (\Fg/\Fh) , M' \right)  &\ra&  
  \Hom_{\FA(\Fh)} (\Cc^{\ot n},  M') \, , \\
\text{resp.} \quad \a^n_{M'}  : \Hom_{\FA(\Fh)} (\Cc^{\ot n},  M')  &\ra&  
\Hom_{\Fh} \left(\bigwedge ^n (\Fg/\Fh) , M' \right) .
\end{eqnarray*}

From (\ref{bM}), one sees that $\, {}_Mb =  \Id \ot  \, b$, where $ b$ 
corresponds to trivial coefficients. In particular, 
\begin{equation*}
{}_Mb \circ {}_M \a \, = \, \Id \ot (b \circ \a) \, = \, 0 \, ,
\end{equation*} 
and likewise, 
\begin{equation*}
{}_M\mu \circ {}_Mb \, = \, \Id \ot (\mu \circ b) \, = \, 0 \, .
\end{equation*} 
Therefore, both
\begin{eqnarray*}
&&{}_M \a : \left( C_* (\Fg, \Fh ; M), 0 \right)  \ra  \left( {\bar C}^*  (\FA(\Fg), \FA(\Fh) ; M) , b \right)
\\
\text{and} \qquad && {}_M \mu  : \left( {\bar C}^*  (\FA(\Fg), \FA(\Fh) ; M) , b \right) \ra
\left( C_* (\Fg, \Fh ; M), 0 \right) \qquad \qquad \qquad \qquad
 \end{eqnarray*} 
 are chain maps with $\, {}_M \mu \circ {}_M \a = \Id$.
 Moreover, as in
the case of trivial coefficients (cf. \cite{Car}, see \cite[XVIII.7]{Ka} for full details)
they induce isomorphism between the corresponding cohomology groups.
\smallskip

On the other hand, for any $\s \in S_{n+1}$ one has
\begin{eqnarray*} 
 \s_{-1}(m \ot_{\Kc} \dot{X}^{\s(1)} \ot \ldots \ot 
 \dot{X}^{\s(n+1)}) = m \dot{X}^{\s(1)} \ot_{\Kc} \dot{X}^{\s(2)} \ot \ldots \ot  \dot{X}^{\s(n+1)}  \\
- \sum_{i=2}^{n+1} m  \ot_{\Kc}  
\dot{X}^{\s(2)} \ot \dots \ot \dot{\widehat{[X^{\s(1)}, X^{\s(i)}]}} \ot \ldots \ot 
 \dot{X}^{\s(n+1)} .
\end{eqnarray*}
By a routine calculation one then sees that, up to a normalizing factor $c_n \neq 0$,
\begin{equation*}
{}_MB \circ {}_M \a \, \simeq \, {}_M\a  \circ  {}_M\d \, .
\end{equation*} 

Thus, when regarding $ \left( C_* (\Fg, \Fh ; M), 0, {}_M\d \right)$ as a bi-complex with
degree $+1$ differential $0$ and degree $-1$ differential ${}_M\d$, $\,  {}_M \a$ is
a bi-chain map, which is a quasi-isomorphism with respect to the degree $+1$ 
differential. From the long exact sequence for $(b, B)$ bi-complexes  \cite{Cext}
it follows that $ {}_M \a$ is also quasi-isomorphism for the corresponding
total complexes.
\end{proof}
\bigskip

Finally, the characteristic map in cyclic cohomology (cf. $\S 4$, Ansatz)
associated an $\Hc$-module algebra $\Ac$ 
with  $\Hc$-invariant trace has a relative version, which we describe below,
in the framework of \cite{HKRS2}. The higher cup products defined
in \cite{KR} should also admit relative counterparts.
\medskip

In addition to the datum of $\S 1$, we consider an $\Hc$-module algebra $\Ac$ 
together with an   \textit{$\Hc$-invariant twisted $M'$-trace} $\phi$ on $\Ac$.
By definition, this means that  
\begin{eqnarray} \label{invt}
\phi &\in& \Hom_{\Hc} (\Ac, M')  
\qquad \qquad \qquad \qquad \qquad \qquad \\ \label{twist} 
\text{and} \qquad \phi( a \, m_{(-1)} (b) )(m_{(0)}) &=& \phi (b \, a) (m),
\qquad a, b \in \Ac , \quad m \in M.
\end{eqnarray}
Examples of such data are provided by actions of Hopf algebras 
with modular pairs $(\d, \s)$ on algebras admitting $\d$-invariant $\s$-traces
and they abound in `nature' (cf. \cite{CMccha, CMcchs}). Furthermore,
the construction of the characteristic map 
associated to such a trace extends
to the present situation in a straightforward fashion. Indeed,
the assignment
\begin{eqnarray*} 
\chi_{\phi}  (m \ot_{\Hc} h^0 \ot \ldots \ot h^n) (a^0,  \ldots , a^n) =
 \phi (h^0(a^0)  h^1 (a^1) \ldots h^n (a^n)) (m)  ,\\ 
  h^0 , \ldots , h^n \in \Hc, \quad
 a^0 , \ldots , a^n \in \Ac , \quad m \in M
\end{eqnarray*}
preserves the cyclic structures and therefore
defines a map from the `absolute' Hopf cyclic cohomology of
$\Hc$ with coefficients in $M$ to the cyclic cohomology of $\Ac$:
\begin{equation}  \label{cmap} 
\chi_{\phi}^* : HC^* (\Hc ; M) \longra HC^* (\Ac) .
\end{equation}
The `relative' version of (\ref{cmap})
arises from the inherited action of the coalgebra
$\Cc = \Cc (\Hc, \Kc)$ on the subalgebra of $\Kc$-invariant elements
$$
\Ac_{\Kc} \, = \, \{ a \in \Ac \, | \quad k(a)\,  =\,  \ve (k) \, a \, , \quad \fl \, k \in \Kc \} \, .
$$
It is defined by a similar formula,
\begin{eqnarray*} 
\chi_{\phi , \Kc}  (m \ot_{\Kc}  c^1 \ot c^2 \ot \ldots \ot c^n) \, = \,
 \phi (a^0  c^1 (a^1) \ldots c^n (a^n)) (m) \, ,\\ \nonumber
 a^0 , \ldots , a^n \in \Ac_{\Kc} , \quad
 c^1 , \ldots , c^n \in \Cc, \quad m \in M ,
\end{eqnarray*}
that again induces at the level of cyclic cohomology the \textit{relative
characteristic map}:
\begin{equation}  \label{rcmap} 
\chi_{\phi, \Kc}^* : HC^* (\Hc , \Kc ; M) \longra HC^* (\Ac_{ \Kc}) .
\end{equation}

\bigskip

\section{The Hopf algebra $\Hc_1$ and its Hopf-cyclic classes}

In this last section we describe in detail the basic Hopf-cyclic 
cocycles of the Hopf algebra $\Hc_1$,
corresponding to the Godbillon-Vey
class, to the Schwarzian derivative and to the transverse
fundamental class,  and illustrate the isomorphism
(\ref{gf}).
\medskip

\noindent We begin by specializing the presentation of the Hopf
algebra $\Hc_n$ (cf. \S 2, Proposition \ref{hopf}) to the 
codimension $1$ case.  
As an algebra, $\Hc_1$ it coincides with
the universal enveloping algebra of the Lie algebra with basis $\{
X,Y,\d_n \, ; n \geq 1 \}$ and brackets
\begin{equation*} \label{pres}
[Y,X] = X \, , \, [Y , \d_n ] = n \, \d_n \, , \, [X,\d_n] =
\d_{n+1} \, , \, [\d_k , \d_{\ell}] = 0 \, , \quad n , k , \ell
\geq 1 \, ,
\end{equation*}
where we have used the abbreviated notation 
$$ X = X_1\, , \, Y = Y^1_{1 1} \, , \, 
\d_1 = \d^1_{1 1} \, , \,  \d_2 = \d^1_{1 1; 1} \, ,
\, \d_3 = \d^1_{1 1; 1 1}  \, , \,  \d_4 = \d^1_{1 1; 1 1 1} \, , \ldots \, .
$$
As a Hopf algebra, the coproduct $\, \D : \Hc_1 \ra \Hc_1 \ot
\Hc_1 \,$ is determined by
\begin{eqnarray*} 
\D \,  Y = Y \ot 1 + 1 \ot Y \, , \quad
\D \,  X &=& X \ot 1 + 1 \ot X + \d_1 \ot Y \nonumber \\
\D \,  \d_1 &=& \d_1 \ot 1 + 1 \ot \d_1
\end{eqnarray*}
and the multiplicativity property
\begin{equation*}
\D (h^1 \, h^2) = \D h^1 \cdot \D h^2 \, , \quad h^1 , h^2 \in
\Hc_1 \, \, ;
\end{equation*}
the antipode is determined by
\begin{equation*}
S(Y) = -Y \, , \, S(X) = -X + \d_1 Y \, , \, S(\d_1) = - \d_1
\end{equation*}
and the anti-isomorphism property
\begin{equation*}
S (h^1 \, h^2) = S(h^2) \, S (h^1) \, , \quad h^1 , h^2 \in \Hc_1
\, \, ;
\end{equation*}
finally, the counit is
\begin{equation*}
\ve (h) = \hbox{constant term of} \quad h \in \Hc_1 \, .
\end{equation*}
\medskip

We next recall the `standard' actions
$\Hc_1$ (cf. \S 2,  Proposition \ref{free}).
Given a one-dimensional manifold $M^1$ and a discrete subgroup $\G
\sbs \Diff^+ (M^1)$, $\, \Hc_1$ acts on the crossed product
algebra
$$\, \Ac_{\G} = C_c^{\ify} (J_+^1 (M^1)) \rtimes \G \, ,
$$
by a Hopf action, where $J_+^1 (M^1) = F^{+} M^1$ is the oriented $1$-jet
bundle over $M^1$. We use the coordinates in $J_+^1 (M^1)$ given
by the Taylor expansion,
$$
j(s) = x + s \, y +  \cdots \, , \qquad y > 0 \, ,
$$
and let diffeomorphisms act in the obvious functorial manner on
the $1$-jets,
$$
{\wt \vp} (x, y) = (\vp (x) , \, \vp' (x) \cdot y) \, .
$$
The canonical action of $\Hc_1$ is then given
 as follows:
\begin{eqnarray*}
Y(f U_{\vp}^*) = y \, \frac{\partial f}{\partial y} \,
U_{\vp}^* \, &,& \qquad X(f U_{\vp}^*) = y \, \frac{\partial
f}{\partial x} \, U_{\vp}^*
\, , \\
\, \d_n (f U_{\vp}^*) &=& y^n \,  \frac{d^{n}}{d x^{n}}
\left(\log \vp' (x) \right) \, f U_{\vp}^* \, \, ,
\label{dn}
\end{eqnarray*}
where we have identified $\,  F^+ M^1 \simeq M^1 \times
\mathbb{R}^+ \, $ and denoted by  $\, (x, y) \, $ the
coordinates on the latter.
\medskip

\noindent The volume form $\displaystyle \frac{dx \wdg
dy}{y^2} \, $ on $ \, F^+ M^1) \,$ is invariant under $\,
\Diff^+ (M^1) \,$ and gives rise to the following trace $\, \tau :
\Ac_{\G} \ra \mathbb{C} $,
\begin{equation} \nonumber
\tau (f U_{\vp}^*) \, = \, \begin{cases} \displaystyle
 \int_{ F^+ M^1)} f(y,y_1) \, \frac{dy \wdg dy_1}{y_1^2}
 &\text{ if }  \vp = 1 \, ,  \\
0 &\text{ if }  \vp \ne 1 \, .
\end{cases}
\end{equation}
This trace is $\d$-invariant with respect to the action $\, \Hc_1
\ot \Ac_{\G} \ra \Ac_{\G} \,$ and with the modular character $\,
\d \in \Hc_1^* \,$, determined by
$$\d (Y) = 1, \quad \d (X) = 0 , \quad \d (\d_n) = 0 \, ;
$$
the invariance property is given by the identity
\begin{equation*}
\tau (h(a)) = \d (h) \, \tau (a) \, , \qquad \fl \quad h \in
\Hc_{1} \, .
\end{equation*}
\smallskip

\noindent The fact that
$$\, S^2 \ne Id \, ,
$$
is automatically corrected by twisting with $\d$. Indeed, $\wt S
= \d * S$ satisfies the involutive property
\begin{equation} \label{inv1}
\wt S^2 = Id \, .
\end{equation}
One has
\begin{equation} \nonumber
\wt S(\delta_1)\, = \, -\delta_1 \,, \quad \wt S(Y) = - Y +1 \, ,
\quad \wt S(X)= -X + \delta_1 Y \, .
\end{equation}
Equation (\ref{inv1}) shows that the pair ($\d,1$) given by the
character $\d$ of $\Hc_1$ and the group-like element $1 \in
\Hc_1$ is a modular pair in involution. Thus
the Hopf-cyclic cohomology $HC^* (\Hc_1; \d, 1)$  is well-defined
and, for each pair $(M^1, \G)$ as above, the assignment
\begin{equation*} 
\chi_{\tau} (h^{1} \ot \ldots \ot h^{n}) (a^0 , \ldots , a^{n}) 
=  \tau (a^0 \, h^{1} (a^{1}) \ldots h^{n} (a^{n})) \, ,
\quad h^i \in \Hc_1 , \, \, a^j \in \Ac \, ,
\end{equation*}
 induces a characteristic homomorphism
\begin{equation} \label{charm} 
 \chi_{\tau}^* \, : \,  HC^* \,(\Hc_1; \d, 1) \, \ra \, HC^* \,
(\Ac_{\G}) \, .
\end{equation}
\medskip

\begin{proposition} \label{gv2}
The element $\, \d_{1} \in \Hc_{1} \,$ is a Hopf cyclic cocycle,
which gives a nontrivial class
$$ 
[\d_{1}] \, \in \, HP^{1} \, (\Hc_1; \d, 1) \, .
$$
\end{proposition}

\begin{proof}
 Indeed, the fact that $\d_{1}$ is a $1$-cocycle is easy to check:
$$
b (\d_1) \, = \, 1 \ot \d_1 - \D \d_1 + \d_1 \ot 1\,  = \, 0 \, ,
$$
while
$$
\tau_1 (\d_1) \, = \, \wt S (\d_1) \, = \, S (\d_1) \, = \, - \d_1
\, .
$$
The non-triviality of to the periodic class $ \, [\d_{1}] \,$ is a consequence
of  Proposition \ref{gv1}  below. Alternatively, one can remark 
 that its image under the above characteristic map (\ref{charm}),
$\, \chi_{\tau}^* \,([\d_{1}]) \, \in \, HC^{1} \, (\Ac_{\G})$, 
is precisely the anabelian $1$-trace of \cite{Ctfc} (cf. also
\cite[III. 6. $\g$]{Cng}), and the latter is known to give a
nontrivial class on the transverse frame bundle to codimension $1$
foliations. 
\end{proof}
\bigskip

\noindent We shall now describe another Hopf cyclic $1$-cocycle,
intimately related to the classical Schwarzian
 \begin{equation} \nonumber
\{y \, ;x\} \, := \, \frac{d^2}{dx^2} \left( \log \frac{dy}{dx}
\right) \, - \, \frac{1}{2} \left( \frac{d}{dx} \left( \log
\frac{dy}{dx} \right) \right)^2 \, ,
\end{equation}
which plays a prominent role in the transverse geometry of
modular Hecke algebras (cf. \cite{CMmha, CMrc}).
\medskip

\begin{proposition} \label{sc1}
The element $\, \d'_{2} := \d_2- \frac{1}{2}\d_1^2 \, \in \Hc_{1}
\,$ is a Hopf cyclic cocycle, whose  action on the crossed product
algebra $ \Ac_{\G} = C_c^{\ify} (J_+^1 (M^1)) \rtimes \G \, $ is
given by the Schwarzian derivative
$$
\d'_2 (f U_{\vp}^*) = y_1^2 \, \{\vp(y) \, ;y\} \, f U_{\vp}^* \, .
$$
Its class
$$  [\d'_{2}] \, \in \, HC^{1} \, (\Hc_1; \d, 1) \,
$$
is equal to $ B[c]$, where $c$ is the following Hochschild
$2$-cocycle:
$$ c \, := \, \delta_1 \ot X + \frac{1}{2} \, \delta_1^2 \ot Y \, .
$$
\end{proposition}
\begin{proof} We shall give the detailed computation, in order
 to illustrate the $(b,B)$ bi-complex for Hopf cyclic cohomology.
Let us compute $b (c)$. One has
\begin{eqnarray*}
b(\delta_1 \ot X) = 1 \ot \delta_1 \ot X - (\delta_1 \ot 1 + 1 \ot
\delta_1) \ot X + \\ \nonumber \delta_1 \ot(X \ot 1 + 1 \ot X +
\delta_1 \ot Y) - \delta_1 \ot X \ot 1  = \delta_1 \ot \delta_1
\ot Y
\end{eqnarray*}
Also
\begin{eqnarray*}
 b(\delta_1^2 \ot Y)= 1 \ot \delta_1^2 \ot Y - (\delta_1^2 \ot 1 +
2 \delta_1 \ot \delta_1 + 1 \ot \delta_1^2) \ot Y + \\ \nonumber
\delta_1^2 \ot ( Y \ot 1 + 1 \ot Y) - \delta_1^2 \ot Y \ot 1 = -2
\delta_1 \ot \delta_1 \ot Y
\end{eqnarray*}
This shows that
$$
    b(c) \, = \, 0 \, ,
$$
so that $c$ is a Hochschild cocycle.
\medskip

\noindent Let us now compute $B(c)$. First, we recall that
$$
B_0 (h^1 \ot h^2) = \wt S (h^1) h^2 \, .
$$
Since $\wt S(\delta_1)\, = \, -\delta_1 \,$, one then has
$$
B_0 (c) = - \delta_1 \, X + \frac{1}{2}  \delta_1^2 \, Y \, .
$$
Since $\, \wt S(Y) = - Y +1 \,$ and $\, \wt S(X)= -X + \delta_1 Y
\,$, it follows that
\begin{eqnarray*}
\wt S(B_0c)  &=&  \wt S(X) \delta_1 + \frac{1}{2}  \wt S(Y)
\delta_1^2  \\ \nonumber &=&(- X + \delta_1 Y)\delta_1 +
\frac{1}{2}  (- Y +1) \delta_1^2  \\ \nonumber
 &=&-X \delta_1 + \delta_1^2 Y + \delta_1^2 - \frac{1}{2}  (\delta_1^2 Y +
\delta_1^2)  \\ \nonumber
  &=&- X \delta_1 + \frac{1}{2}  \delta_1^2 Y + \frac{1}{2}  \delta_1^2 \, .
\end{eqnarray*}
Therefore,
\begin{eqnarray*}
B(c)  \, = \, B_0c - \wt S(B_0c)  \, &=& \,
\\ \nonumber
- \delta_1 \, X + \frac{1}{2}  \delta_1^2 \, Y - (-X \delta_1 +
\frac{1}{2}  \delta_1^2 Y + \frac{1}{2}  \delta_1^2) \,
 \, &=& \, \delta'_2 \, .
\end{eqnarray*}
which shows that the class of $\d'_2$ is trivial in the periodic
cyclic cohomology $HP^{*}(\Hc_1; \d, 1)$.
\end{proof}
\medskip

Let us illustrate the isomorphism (\ref{Tnat})  
 by determining the $2$-cocycle that represents the `universal' transversal
fundamental class in codimension $1$.
\medskip

\begin{proposition} \label{utfc}
 The  $2$-cochain
\begin{equation} \nonumber
\Pi := X \ot Y - Y \ot X - \d_1 \,  Y \ot Y \, \in \, \Hc_1 \ot \Hc_1
\end{equation}
is a cyclic $2$-cocycle, whose class $[\Pi] \in HC^{2}(\Hc_1; \d, 1)$
 corresponds by $T^{\natural}$
to the transverse fundamental class.
\end{proposition}
\smallskip

\begin{proof} The transverse fundamental class in $HC^2 (\Ac_1)$
is given by the $2$-cocycle
\begin{equation} \nonumber
F(f^0 U^*_{\vp_0} , f^1 U^*_{\vp_1} , f^2 U^*_{\vp_2}) =  
\left\{ \begin{matrix} \displaystyle
\int_{F^{+} \Rb} f^0 \, {\wt \vp_0}^* ( d f^1) \,  {\wt \vp_0}^* {\wt \vp_1}^* (d f^2)
\, , \quad \text{if}
\quad \vp_2  \vp_1 \vp_0 = Id \, , \cr\cr
\quad 0 \, , \qquad \qquad  \qquad \qquad \text{otherwise} \, .
\end{matrix} \right.
\end{equation}
In terms of the standard basis of $T^*(F^{+} \Rb)$, one has
\begin{equation} \nonumber
df^1 \, = \, X(f^1) \, y^{-1} dx \, + \, Y(f^1)\,  y^{-1} dy  \, ,
\end{equation}
 therefore
\begin{eqnarray} \label{phif} \nonumber
 {\wt \vp_0}^* ( d f^1) &=& {\wt \vp_0}^*\left(X(f^1)\right) \cdot  {\wt \vp_0}^* (y^{-1} dx) \, + \, 
 {\wt \vp_0}^*\left(Y(f^1)\right)\cdot  {\wt \vp_0}^* ( y^{-1} dy )  \\ \nonumber
&=& \left( {\wt \vp_0}^*(X(f^1)) + \g_1(\vp_0) \,  {\wt \vp_0}^*(Y(f^1))\right)\,  y^{-1} dx  +
 {\wt \vp_0}^*\left(Y(f^1)\right)\, y^{-1} dy  \, ,
\end{eqnarray}
where, consistently with the abbreviation $\d_1$,  $\, \g_1 = \g^1_{1 1}$.

On substituting the above expression of $ {\wt \vp_0}^* ( d f^1)$ together with the similar one
for $ {\wt \vp_0}^* {\wt \vp_1}^* (d f^2)$ in the formula defining $F$, after using
the cocycle identity
$$
\g_1(\vp_1 \vp_0) \, = \, {\wt \vp_0}^* ( \g_1 (\vp_1)) + \g_1 (\vp_0) \, ,
$$
one observes that the result is identical to
$\, T^{\natural} (\Pi) (f^0 U^*_{\vp_0} , f^1 U^*_{\vp_1} , f^2 U^*_{\vp_2}) $.
\end{proof}
\bigskip

We now proceed to
illustrate the isomorphism
\begin{equation} \label{GF1}
\kappa_{1}^{*} \, : \sum^{\oplus}_{i \equiv * \, (2)} \quad  H^{i} ({\Fa}_1 , \mathbb{C})
\build\longrightarrow_{}^{\simeq} HP^{*} \, (\Hc_1 ; \d, 1)
\end{equation}
between the Gelfand-Fuks cohomology of the Lie algebra of formal
vector fields on $\mathbb{R}^{1}$ and the periodic Hopf-cyclic
cohomology of $(\Hc_1; \d, 1)$.

We recall (cf. \cite{GF, Go}) that the cohomology $\, H^{*} ({\Fa}_1 , \mathbb{R}) \, $
is finite dimensional and the only nontrivial groups are:
$$
H^0 ({\Fa}_1 , \mathbb{R}) = \mathbb{R} \cdot 1 \qquad \hbox{and}
\qquad H^3 ({\Fa}_1 , \mathbb{R}) = \mathbb{R} \cdot gv \, ,
$$
where
\begin{equation*}
gv (p_1 \, \partial_x , p_2 \, \partial_x , p_3 \, \partial_{x}) =
\left\vert \begin{array}{ccc}
p_1 (0) &p_2 (0) &p_3 (0) \\ \\
p'_1 (0) &p'_2 (0) &p'_3 (0) \\ \\
p''_1 (0) &p''_2 (0) &p''_3 (0)
\end{array} \right\vert \, .
\end{equation*}
\medskip

\noindent   The
class of the unit constant $\, [1] \in HC^{0}(\Hc_1; \d, 1) \,$
is  trivial in $HP^{0}(\Hc_1; \d, 1) \,$, since $ \, B (Y) = 1 \,$.
On the other hand,
\begin{equation} \label{fund}
\kappa_{1}^{*} (1) \, = \, F \, ,
\end{equation}
which thus gives the generator of
 $HP^{0}(\Hc_1; \d, 1)$. Formula (\ref{fund}) is easy to check
 using Proposition \ref{utfc} and applying the
 definition of $\kappa_{1}^{*}$, cf. \cite{CMhti} or \cite{CMdcc},  
  in the case of $0$-dimensional Lie algebra cochains.
\medskip

On the other hand, the evaluation of $\kappa_{1}^{*}$ on the
the Godbillon-Vey class requires some calculations, whose
details will be given below, after we recall the original definition
of this class.
\smallskip

\noindent Let $V$ be a closed, smooth manifold, foliated by a
transversely oriented codimension $1$ foliation $\Fc$. Then $T\Fc
= \Ker \om \sbs TV$, for some $\om \in \Om^1 (V)$ such that $\om
\wdg d\om = 0$. Equivalently, $d\om = \om \wdg \a$ for some $\a
\in \Om^1 (V)$, which implies $d\a \wdg \om = 0$. In turn, the
latter ensures that $d\a = \om \wdg \b$, $\b \in \Om^1 (V)$. Thus,
$\a \wdg d\a \in \Om^3 (V)$ is closed. Its de~Rham cohomology
class,
\begin{equation*}
\GV (V, \Fc) = [\a \wedge d\a] \in H^3 (V,\mathbb{R}) \, ,
\end{equation*}
is independent of the choices of $\om$ and $\a$ and represents the
original definition of the Godbillon-Vey class.
\smallskip

\noindent The Godbillon-Vey class acquires a universal status when
viewed as a characteristic class (cf. \cite{H1}) associated to the
Gelfand-Fuks cohomology of the Lie algebra ${\Fa}_1 =
\mathbb{R} \, [[x]] \, \partial_x$ of formal vector fields on
$\mathbb{R}$. 
\smallskip

\noindent Given any oriented $1$-dimensional manifold $M^1$, the
Lie algebra cocycle $gv$ can be converted into a $3$-form on the
jet bundle (of orientation preserving jets)
$$ J_+^{\ify} (M^1) = \limproj_n J_+^n (M^1) \, ,
$$
invariant under the pseudogroup $\Gc^+ (M^1)$ of all orientation
preserving local diffeomorphisms of $M^1$. Indeed, sending the
formal vector field
$$
p = j_0^{\ify} \left( \frac{dh_t}{dt} \Biggl\vert_{t=0} \right)
\in \Fa_1 \, ,
$$
where $\{ h_t \}$ is a $1$-parameter family of local
diffeomorphisms of $\mathbb{R}$ preserving the origin, to the
$\Gc^+ (M^1)$-invariant vector field
$$
j_0^{\ify} \left( \frac{d(f \circ h_t)}{dt} \Biggl\vert_{t=0}
\right) \in T_{j_0^{\ify} (f)} \, J_+^{\ify} (M^1) \,
$$
gives a natural identification of the Lie algebra complex of
${\Fa}_1$ with the invariant forms on the jet bundle,
$$
 \quad \t \in C^{\bullet} ({\Fa}_1) \mapsto \t \in \Omega^{\bullet}
(J_+^{\ify} (M^1))^{\Gc^+ (M^1)} \, .
$$
In local coordinates on $J_+^{\ify} (M^1)$, given by the
coefficients of the Taylor expansion at $0$,
$$
f(s) = x + s \, y + s^2 y_1 + \cdots \, , \qquad y > 0 \, ,
$$
one has
\begin{eqnarray}
dx \, &=& \, y \, \t^0  \nonumber \\
dy \, &=& \, y \, \t^1 + 2 \, y_1 \, \t^0
                        \nonumber \\
dy_1 \, &=& \, y \, \t^2 + 2 \, y_1 \, \t^1 + 3 \, y_2 \, \t^0
\, , \nonumber
\end{eqnarray}
therefore
\begin{equation*}
 \quad gv \equiv \t^0 \wdg \t^1 \wdg \t^2 \, = \,  \frac{1}{y^3} \, dx \wdg dy  \wdg dy_1 \in
\Om^3 (J_+^{\ify} (M^1))^{\Gc^+ (M^1)} \, .
\end{equation*}
\medskip

\noindent Given a codimension $1$ foliation $(V, \Fc)$ as above,
one can find an open covering $\{ U_i \}$ of $V$ and submersions
$f_i : U_i \ra T_i \sbs \mathbb{R}$, whose fibers are plaques of
$\Fc$, such that $f_i = g_{ij} \circ f_j$ on $U_i \cap U_j$, with
$g_{ij} \in \Gc^+ (M^1)$ a $1$-cocycle. Then $\, M^1 = \bigcup_i
T_i \times \{i\} \, $ is a complete transversal. Let $J^{\infty}
(\Fc)$ denote the bundle over $V$ whose fiber at $x \in U_i$
consists of the $\infty$-jets of local submersions of the form
$\varphi \circ f_i$ with $\varphi \in \Gc^+ (M^1)$. Using the
$\Gc^+ (M^1)$-invariance of $gv \in \Om^3 (J^{\infty} (M^1))$ one
can pull it back to a closed form $gv (\Fc) \in \Om^3 (J^2
(\Fc))$. Its de~Rham class  $ [gv (\Fc)] \in H^3 (J^2 (\Fc),
\mathbb{R})$, when viewed as a class in $H^3 (V, \mathbb{R})$ (the
fibers of $J^2 \Fc$ being contractible), is precisely the
Godbillon-Vey class $\GV (V, \Fc)$.

\begin{proposition} \label{gv1} The canonical cochain map
associated to the trivial connection on $\, F^+ {\Rb} \,$
sends the universal Godbillon-Vey cocycle $\, gv \,$  to the Hopf
cyclic cocycle $ \, \d_{1} \,$, implementing the identity
\begin{equation*}
\kappa_{1}^{*} ([gv]) \, = \, [\d_{1}] \, .
\end{equation*}
\end{proposition}
\smallskip

\begin{proof} The definition of $\kappa_{1}^{*}$ involves two steps (see
\cite{CMhti}, \cite{CMdcc} ). The first turns
the Lie algebra
cocycle $gv \in C^{2} ({\Fa}_1, \mathbb{R})$,  or equivalently,
the  form on the jet bundle $J_+^2 (\mathbb{R})$
$$
gv \, = \, \frac{1}{y^3} \, dx \wdg dy \wdg dy_1 \, ,
$$
 into a group
$1$-cocycle $\, C_{1,0} (gv) \, $ on $\, \Gc_1 = \Diff^+
(\mathbb{R}) \,$ with values in currents on $\, J_+^1 (\mathbb{R})
\,$.  The second step takes its image in the cyclic bi-complex
 under the canonical map $\Phi$. One thus obtains a cyclic
cocycle
$$
(\Phi \, ( C_{1,0} \, (gv) \,) ) \, (f^0 \, U_{\vp_0}^* , f^1 \,
U_{\vp_1}^*)
$$
which is automatically supported at the identity, i.e. it is nonzero
only when $\vp_1 \vp_0 = 1$. Moreover, it is of the form
$$
(\Phi \, ( C_{1,0} \, (gv) \, ) ) \, (f^0 \, U_{\vp}^* , f^1 \,
U_{\vp^{-1}}^*) \, = \, -
  \langle \, C_{1,0} (gv) (1, \vp) , \, f^0 \cdot {\wt \vp}^* (f^1) \,  \rangle \, .
$$
Applying the definition given in \cite{CMdcc}, one gets
\begin{equation*}
\langle C_{1,0} (gv) (1 , \vp) , f \rangle \, = \, \int_{\D^1 \ts
F^+ {\Rb}} f \, \wt{\s} (1 , \vp)^* \left( gv \right)
\end{equation*}
where $\, \D^1 \,$ is the $1$-simplex and
$$ \wt{\s} (1 , \vp)  : \D^1 \ts J_+^1 (\mathbb{R} )
  \ra J_+^{\ify} (\mathbb{R})
$$
has the expression
$$
\wt{\s} (1 , \vp) (t , x , y) \, = \, \s_{(1-t) \nb_{0} + t
\nb_{0}^{\vp}} (x , y) \, ,
$$
where the meaning of the notation used is as follows.

First, $\, \nb_{0} \,$ stands for the trivial linear connection on
$\mathbb{R}$, given by the connection form on  $F^+ {\Rb}$
$$
\om_{0} = y^{-1} \, dy \, ,
$$
while $\, \nb_{0}^{\vp} \,$  denotes its transform under the
prolongation
$$
{\wt \vp} (x, y) = (\vp (x) , \, \vp' (x) \cdot y) \, ,
$$
of the diffeomorphism $ \vp \in  \Gc_1 $\, ; the latter corresponds
to the connection form
$$
{\wt \vp}^* (\om_{0}) = \frac{1}{\vp' (x) \, y} \, \left(\vp' (x) \, dy +
\vp'' (x) \cdot y \, dx\right) = y^{-1} \, dy + \frac{d}{dx} (\log \vp' (x)) \, dx \, .
$$
Furthermore, for any linear connection $\, \nb \,$ on
$\mathbb{R}$, $ \, \s_{\nb} \,$ denotes the jet
$$
\s_{\nb} (x, y) = j_0^{\ify} \left( Y (s) \right) \,
$$
of the local diffeomorphism
$$ s \, \mpo  \, Y(s) := \exp_y^{\nb} \left( s \, y \, \frac{d}{dx}
\right) \, , \qquad s \in \mathbb{R} \, .
$$
Now $ \, Y(s) \,$ satisfies the geodesics ODE
$$
\left\{ \begin{matrix} \ddot Y (s) + \G_{1 \,1}^{1} (Y(s)) \cdot
\dot Y (s)^2 = 0 \, , \hfill \cr Y(0) = x \, , \hfill \cr \dot Y
(0) = y \, . \hfill \cr
\end{matrix} \right. \,
$$
Since we only need the 2-jet of the exponential map, suffices to
retain that
$$
Y (0) = x \, , \quad \dot Y (0) = y \quad \hbox{and} \quad \ddot
Y (0) = - \G_{1 \,1}^{1} (x) \, y^2 \, .
$$
Thus,
$$
\s_{\nb} (x, y) \, = \, x + y \, s \, - \, \G_{1 \,1}^{1} (x)
\, y^2 \, s^2  \, + \, \hbox{higher order terms} \, .
$$

In our case  $\, \nb \, = \, (1-t) \nb_{0} \,  + \, t
\nb_{0}^{\vp} \, $, which gives
$$
\G_{1 \,1}^{1} (t, x) \, = \, t \, \frac{d}{dx} ( \log \vp' (x))
\, ,
$$
and therefore
$$
\wt{\s} (1 , \vp) (t , x , y) \, = \, \, x + y \, s \, - \, t
\, \frac{d}{dx} ( \log \vp' (x))
 \, y^2 \, s^2  \, + \, \hbox{higher order terms} \, .
$$
It follows that
$$
\wt{\s} (1 , \vp)^* (dx) \, = \, dx \, , \qquad \wt{\s} (1 ,
\vp)^* (dy) \, = \, dy
$$
and
\begin{eqnarray}
\wt{\s} (1 , \vp)^* (dy_1) \, = \, - \left(  \frac{d}{dx} ( \log
\vp' (x)) \, dt \, + \, t \frac{d}{d x} \left(\frac{d}{dx} ( \log
\vp' (x))\right)
\, dx \right) y^2  \nonumber \\
 - \, 2 t \, \frac{d}{dx} ( \log \vp' (x)) \, y
\, d y \, . \nonumber
\end{eqnarray}
Hence on $\D^1 \ts F^+ {\Rb}$,
$$
\wt{\s} (1 , \vp)^*  (gv) \,  = \, - \frac{1}{y} \, \frac{d}{dx}
 ( \log \vp' (x))
\, dt \wdg dx \wdg dy \, .
$$
Going back to the definition of the group cochain, one gets
\begin{eqnarray*}
\langle C_{1,0} (gv) (1 , \vp) , f \rangle \, &=& - \int_{ F^+ {\Rb}}
 f (x , y) \cdot \int_0^1
 dt \cdot \frac{1}{y} \, \frac{d}{dx} ( \log \vp' (x))
 dx \wedge dy  \\
&= & - \int_{ F^+ {\Rb} }  f(x , y) \left( y \
\frac{d}{dx} ( \log \vp' (x)) \right) \, \frac{dx \wedge
dy}{y^2}   \, ,
\end{eqnarray*}
which finally gives
\begin{eqnarray}
&\ &(\, \Phi \, ( C_{1,0} \, (gv)))
(f^0 \, U_{\vp}^* , f^1 \, U_{\vp^{-1}}^*) \, = \nonumber \\
&=& \, \int_{ F^+ {\Rb}} f^0 \cdot {\wt \vp}^* (f^1)  \cdot
\left( y \, \frac{d}{dx} ( \log \vp' (x)) \right) \,
\frac{dx \wedge dy}{y^{2}} \nonumber \\
&=& \tau (f^0 U_{\vp}^* \cdot \d_1 (f^1 U_{\vp^{-1}}^*)) \, = \, \,
\chi_{\tau} (\d_1) (f^0 \, U_{\vp}^* , f^1 \, U_{\vp^{-1}}^*) \, .
\nonumber
\end{eqnarray}
\end{proof}

 \end{document}